% formatam.tex -- AMSTeX template file
% Version: June 10, 1992 

\input amstex

%\documentstyle{amams} % input Annals of Mathematics macros.
%%%%%%%%%%%%%%%%%%%%%%%%%%%%%%%%%%%%%%%%%%%%%%%%%%%%%%%%%%%
%% amams.sty: AMSTeX Macros for Articles to be published in
%%
%% Annals of Mathematics
%% 
%% Princeton University and the
%% Institute for Advanced Study
%%
%% Published by Princeton University Press
%%%%%%%%%%%%%%%%%%%%%%%%%%%%%%%%%%%%%%%%%%%%%%%%%%%%%%%%%%%

%%%%%%%%%%%%%%%%%%%%%%%%%%%%%%%%%%%%%%%%%%%%%%%%%%%%%%%%%%%
%% Variations on AMSPPT.sty written by Amy Hendrickson
%% TeXnology Inc, Brookline, MA
%% 617 738-8029, amyh@ai.mit.edu
%%%%%%%%%%%%%%%%%%%%%%%%%%%%%%%%%%%%%%%%%%%%%%%%%%%%%%%%%%%

\def\spaces{\space\space\space\space\space\space\space\space\space\space}
\def\spacess{\message{\spaces\spaces\spaces\spaces\spaces\spaces\spaces}}
\spacess
\spacess
\message{Annals of Mathematics Style: Current Version: 1.1. June 10, 1992}
\spacess
\spacess
%%%%%%%%%%%%%%%%%%%%%%%%%%%%%%%%%%%%%%%%%%%%%%%%%%%%%%%%

\catcode`\@=11

\hyphenation{acad-e-my acad-e-mies af-ter-thought anom-aly anom-alies
an-ti-deriv-a-tive an-tin-o-my an-tin-o-mies apoth-e-o-ses apoth-e-o-sis
ap-pen-dix ar-che-typ-al as-sign-a-ble as-sist-ant-ship as-ymp-tot-ic
asyn-chro-nous at-trib-uted at-trib-ut-able bank-rupt bank-rupt-cy
bi-dif-fer-en-tial blue-print busier busiest cat-a-stroph-ic
cat-a-stroph-i-cally con-gress cross-hatched data-base de-fin-i-tive
de-riv-a-tive dis-trib-ute dri-ver dri-vers eco-nom-ics econ-o-mist
elit-ist equi-vari-ant ex-quis-ite ex-tra-or-di-nary flow-chart
for-mi-da-ble forth-right friv-o-lous ge-o-des-ic ge-o-det-ic geo-met-ric
griev-ance griev-ous griev-ous-ly hexa-dec-i-mal ho-lo-no-my ho-mo-thetic
ideals idio-syn-crasy in-fin-ite-ly in-fin-i-tes-i-mal ir-rev-o-ca-ble
key-stroke lam-en-ta-ble light-weight mal-a-prop-ism man-u-script
mar-gin-al meta-bol-ic me-tab-o-lism meta-lan-guage me-trop-o-lis
met-ro-pol-i-tan mi-nut-est mol-e-cule mono-chrome mono-pole mo-nop-oly
mono-spline mo-not-o-nous mul-ti-fac-eted mul-ti-plic-able non-euclid-ean
non-iso-mor-phic non-smooth par-a-digm par-a-bol-ic pa-rab-o-loid
pa-ram-e-trize para-mount pen-ta-gon phe-nom-e-non post-script pre-am-ble
pro-ce-dur-al pro-hib-i-tive pro-hib-i-tive-ly pseu-do-dif-fer-en-tial
pseu-do-fi-nite pseu-do-nym qua-drat-ics quad-ra-ture qua-si-smooth
qua-si-sta-tion-ary qua-si-tri-an-gu-lar quin-tes-sence quin-tes-sen-tial
re-arrange-ment rec-tan-gle ret-ri-bu-tion retro-fit retro-fit-ted
right-eous right-eous-ness ro-bot ro-bot-ics sched-ul-ing se-mes-ter
semi-def-i-nite semi-ho-mo-thet-ic set-up se-vere-ly side-step sov-er-eign
spe-cious spher-oid spher-oid-al star-tling star-tling-ly
sta-tis-tics sto-chas-tic straight-est strange-ness strat-a-gem strong-hold
sum-ma-ble symp-to-matic syn-chro-nous topo-graph-i-cal tra-vers-a-ble
tra-ver-sal tra-ver-sals treach-ery turn-around un-at-tached un-err-ing-ly
white-space wide-spread wing-spread wretch-ed wretch-ed-ly Brown-ian
Eng-lish Euler-ian Feb-ru-ary Gauss-ian Grothen-dieck Hamil-ton-ian
Her-mit-ian Jan-u-ary Japan-ese Kor-te-weg Le-gendre Lip-schitz
Lip-schitz-ian Mar-kov-ian Noe-ther-ian No-vem-ber Rie-mann-ian
Schwarz-schild Sep-tem-ber Za-mo-lod-chi-kov Knizh-nik quan-tum Op-dam
Mac-do-nald Ca-lo-ge-ro Su-ther-land Mo-ser Ol-sha-net-sky  Pe-re-lo-mov
in-de-pen-dent ope-ra-tors Hec-man Op-dam har-mo-nic ana-ly-sis
al-geb-ra pre-print char-ac-ter-is-tic sub-mo-dule re-la-tion
}

\Invalid@\nofrills
\Invalid@\usualspace
\newif\ifnofrills@
\def\nofrills@#1#2{\relaxnext@
  \DN@{\ifx\next\nofrills
    \nofrills@true\let#2\relax\DN@\nofrills{\nextii@}%
  \else
    \nofrills@false\def#2{#1}\let\next@\nextii@\fi
\next@}}
\def\usualspace@#1{\ifnofrills@\def\usualspace{#1}\fi}
\def\addto#1#2{\csname \expandafter\eat@\string#1@\endcsname
  \expandafter{\the\csname \expandafter\eat@\string#1@\endcsname#2}}
\newdimen\bigsize@
\def\big@#1#2{{\hbox{$\left#2\vcenter to#1\bigsize@{}%
  \right.\nulldelimiterspace\z@\m@th$}}}
\def\big{\big@\@ne}
\def\Big{\big@{1.5}}
\def\bigg{\big@\tw@}
\def\Bigg{\big@{2.5}}
\def\raggedcenter@{\leftskip\z@ plus.4\hsize \rightskip\leftskip
 \parfillskip\z@ \parindent\z@ \spaceskip.3333em \xspaceskip.5em
 \pretolerance9999\tolerance9999 \exhyphenpenalty\@M
 \hyphenpenalty\@M \let\\\linebreak}
\def\upperspecialchars{\def\ss{SS}\let\i=I\let\j=J\let\ae\AE\let\oe\OE
  \let\o\O\let\aa\AA\let\l\L}
\def\uppercasetext@#1{%
  {\spaceskip1.2\fontdimen2\the\font plus1.2\fontdimen3\the\font
   \upperspecialchars\uctext@#1$\m@th\aftergroup\eat@$}}
\def\uctext@#1$#2${\endash@#1-\endash@$#2$\uctext@}
\def\endash@#1-#2\endash@{%
\uppercase{#1}\if\notempty{#2}--\endash@#2\endash@\fi}
\def\runaway@#1{\DN@{#1}\ifx\envir@\next@
  \Err@{You seem to have a missing or misspelled \string\end#1 ...}%
  \let\envir@\empty\fi}
\newif\iftemp@
\def\notempty#1{TT\fi\def\test@{#1}\ifx\test@\empty\temp@false
  \else\temp@true\fi \iftemp@}

%\comment%%% remove
\font@\tensmc=cmcsc10
\font@\sevenex=cmex7 
\font@\sevenit=cmti7
\font@\eightrm=cmr8 % preloaded in plain.tex
\font@\sixrm=cmr6 % preloaded in plain.tex
\font@\eighti=cmmi8     \skewchar\eighti='177 % preloaded
\font@\sixi=cmmi6       \skewchar\sixi='177   % preloaded
\font@\eightsy=cmsy8    \skewchar\eightsy='60 % preloaded
\font@\sixsy=cmsy6      \skewchar\sixsy='60   % preloaded
\font@\eightex=cmex8 %
\font@\eightbf=cmbx8 % preloaded in plain.tex
\font@\sixbf=cmbx6   % preloaded in plain.tex
\font@\eightit=cmti8 % preloaded in plain.tex
\font@\eightsl=cmsl8 % preloaded in plain.tex
\font@\eightsmc=cmcsc10 
\font@\eighttt=cmtt8 % preloaded in plain.tex
%\font@\ninerm=cmr9
%\font@\ninei=cmmi9    \skewchar\ninei='177
%\font@\ninesy=cmsy9   \skewchar\ninesy='60
%\font@\nineex=cmex9
%\font@\ninebf=cmbx9
%\font@\nineit=cmti9
%\font@\ninesl=cmsl9
%\font@\ninesmc=cmcsc9
%\font@\ninemsa=msam9
%\font@\ninemsb=msbm9
%\font@\nineeufm=eufm9
%\endcomment%%%

\loadmsam
\loadmsbm
\loadeufm
\UseAMSsymbols

\def\penaltyandskip@#1#2{\relax\ifdim\lastskip<#2\relax\removelastskip
      \ifnum#1=\z@\else\penalty@#1\relax\fi\vskip#2%
  \else\ifnum#1=\z@\else\penalty@#1\relax\fi\fi}
\def\nobreak{\penalty\@M
  \ifvmode\def\penalty@{\let\penalty@\penalty\count@@@}%
  \everypar{\let\penalty@\penalty\everypar{}}\fi}
\let\penalty@\penalty

\def\block{\RIfMIfI@\nondmatherr@\block\fi
       \else\ifvmode\vskip\abovedisplayskip\noindent\fi
        $$\def\endblock{\par\egroup$$}\fi
  \vbox\bgroup\advance\hsize-2\indenti\noindent}
\def\endblock{\par\egroup}

\def\logo@{\baselineskip2pc \hbox to\hsize{\hfil\eightpoint Typeset by
 \AmSTeX}}

%%%%%%%%%%%%%%%%%%%%%%%%%%%%%%%%%%%%%%%%%%%%%%%%%%%%%%%%%%%%%%%
%% Macros for Annals of Mathematics written by Amy Hendrickson
%% TeXnology Inc, Brookline, MA
%% 617 738-8029, amyh@ai.mit.edu
%%%%%%%%%%%%%%%%%%%%%%%%%%%%%%%%%%%%%%%%%%%%%%%%%%%%%%%%%%%%%%%

%% This file includes: 
%% 1) Font declarations, 
%% 2) Page set up, 
%% 3) Title page 
%% 4) Section heads,
%% 5) Equation macros, autonumbering equations, etc.,
%% 6) Figure and Table Captions,
%% 7) End matter macros: Bibliography, Appendix, etc.,
%% 8) Footnotes,
%% 9) Theorem type environments
%% 10) Cross-referencing
%% 11) Listing
%% 12) Article and Journal Table of Contents

%%%%%%%%%%%%%%%%%%%%%%%%%%%%%%%%%%%
%% 1) Font declarations, 
% Computer Modern fonts 

% Small Caps
\font\elevensc=cmcsc10 scaled\magstephalf
\font\tensc=cmcsc10

\font\eightsc=cmcsc10 scaled800

\font\elevenrm=cmr10 scaled \magstephalf%!!!
\font\ninerm=cmr9
\font\eightrm=cmr8
\font\sixrm=cmr6
\font\fiverm=cmr5

\font\eleveni=cmmi10 scaled\magstephalf
\font\ninei=cmmi9
\font\eighti=cmmi8
\font\sixi=cmmi6
\font\fivei=cmmi5
\skewchar\ninei='177 \skewchar\eighti='177 \skewchar\sixi='177
\skewchar\eleveni='177

\font\elevensy=cmsy10 scaled\magstephalf
\font\ninesy=cmsy9
\font\eightsy=cmsy8
\font\sixsy=cmsy6
\font\fivesy=cmsy5
\skewchar\ninesy='60 \skewchar\eightsy='60 \skewchar\sixsy='60
\skewchar\elevensy'60

\font\eighteenbf=cmbx10 scaled\magstep3

\font\twelvebf=cmbx10 scaled \magstep1
\font\elevenbf=cmbx10 scaled \magstephalf
\font\tenbf=cmbx10
\font\ninebf=cmbx9
\font\eightbf=cmbx8
\font\sixbf=cmbx6
\font\fivebf=cmbx5

\font\elevenit=cmti10 scaled\magstephalf
\font\nineit=cmti9
\font\eightit=cmti8

% Fonts for bold math
\font\eighteenmib=cmmib10 scaled \magstep3
\font\twelvemib=cmmib10 scaled \magstep1
\font\elevenmib=cmmib10 scaled\magstephalf
\font\tenmib=cmmib10
\font\eightmib=cmmib10 scaled 800 
\font\sixmib=cmmib10 scaled 600

\font\eighteensyb=cmbsy10 scaled \magstep3 
\font\twelvesyb=cmbsy10 scaled \magstep1
\font\elevensyb=cmbsy10 scaled \magstephalf
\font\tensyb=cmbsy10 
\font\eightsyb=cmbsy10 scaled 800
\font\sixsyb=cmbsy10 scaled 600
 
\font\elevenex=cmex10 scaled \magstephalf
\font\tenex=cmex10     
\font\eighteenex=cmex10 scaled \magstep3

%%%%%%%%%%%%%%%%%%%%%%%%%%%%
%% Font families

\def\elevenpoint{\def\rm{\fam0\elevenrm}%
  \textfont0=\elevenrm \scriptfont0=\eightrm \scriptscriptfont0=\sixrm
  \textfont1=\eleveni \scriptfont1=\eighti \scriptscriptfont1=\sixi
  \textfont2=\elevensy \scriptfont2=\eightsy \scriptscriptfont2=\sixsy
  \textfont3=\elevenex \scriptfont3=\tenex \scriptscriptfont3=\tenex
  \def\bf{\fam\bffam\elevenbf}%
  \def\it{\fam\itfam\elevenit}%
  \textfont\bffam=\elevenbf \scriptfont\bffam=\eightbf
   \scriptscriptfont\bffam=\sixbf
\normalbaselineskip=13.95pt
  \setbox\strutbox=\hbox{\vrule height9.5pt depth4.4pt width0pt\relax}%
  \normalbaselines\rm}

\elevenpoint %%% default fonts and baselineskip

\def\ninepoint{\def\rm{\fam0\ninerm}%
  \textfont0=\ninerm \scriptfont0=\sixrm \scriptscriptfont0=\fiverm
  \textfont1=\ninei \scriptfont1=\sixi \scriptscriptfont1=\fivei
  \textfont2=\ninesy \scriptfont2=\sixsy \scriptscriptfont2=\fivesy
  \textfont3=\tenex \scriptfont3=\tenex \scriptscriptfont3=\tenex
  \def\it{\fam\itfam\nineit}%
  \textfont\itfam=\nineit
  \def\bf{\fam\bffam\ninebf}%
  \textfont\bffam=\ninebf \scriptfont\bffam=\sixbf
   \scriptscriptfont\bffam=\fivebf
\normalbaselineskip=11pt
  \setbox\strutbox=\hbox{\vrule height8pt depth3pt width0pt\relax}%
  \normalbaselines\rm}

\def\eightpoint{\def\rm{\fam0\eightrm}%
  \textfont0=\eightrm \scriptfont0=\sixrm \scriptscriptfont0=\fiverm
  \textfont1=\eighti \scriptfont1=\sixi \scriptscriptfont1=\fivei
  \textfont2=\eightsy \scriptfont2=\sixsy \scriptscriptfont2=\fivesy
  \textfont3=\tenex \scriptfont3=\tenex \scriptscriptfont3=\tenex
  \def\it{\fam\itfam\eightit}%
  \textfont\itfam=\eightit
  \def\bf{\fam\bffam\eightbf}%
  \textfont\bffam=\eightbf \scriptfont\bffam=\sixbf
   \scriptscriptfont\bffam=\fivebf
\normalbaselineskip=12pt
  \setbox\strutbox=\hbox{\vrule height8.5pt depth3.5pt width0pt\relax}%
  \normalbaselines\rm}

%%%%%%%%%%%%%%%%%%%%%%%%%%%%
%% Font families for bold math in title and section heads

\def\eighteenbold{\def\rm{\fam0\eighteenbf}%
  \textfont0=\eighteenbf \scriptfont0=\twelvebf \scriptscriptfont0=\tenbf
  \textfont1=\eighteenmib \scriptfont1=\twelvemib\scriptscriptfont1=\tenmib
  \textfont2=\eighteensyb \scriptfont2=\twelvesyb\scriptscriptfont2=\tensyb
  \textfont3=\eighteenex \scriptfont3=\tenex \scriptscriptfont3=\tenex
  \def\bf{\fam\bffam\eighteenbf}%
  \textfont\bffam=\eighteenbf \scriptfont\bffam=\twelvebf
   \scriptscriptfont\bffam=\tenbf
\normalbaselineskip=20pt
  \setbox\strutbox=\hbox{\vrule height13.5pt depth6.5pt width0pt\relax}%
\everymath {\fam0 }
\everydisplay {\fam0 }
  \normalbaselines\rm}

\def\elevenbold{\def\rm{\fam0\elevenbf}%
  \textfont0=\elevenbf \scriptfont0=\eightbf \scriptscriptfont0=\sixbf
  \textfont1=\elevenmib \scriptfont1=\eightmib \scriptscriptfont1=\sixmib
  \textfont2=\elevensyb \scriptfont2=\eightsyb \scriptscriptfont2=\sixsyb
  \textfont3=\elevenex \scriptfont3=\elevenex \scriptscriptfont3=\elevenex
  \def\bf{\fam\bffam\elevenbf}%
  \textfont\bffam=\elevenbf \scriptfont\bffam=\eightbf
   \scriptscriptfont\bffam=\sixbf
\normalbaselineskip=14pt
  \setbox\strutbox=\hbox{\vrule height10pt depth4pt width0pt\relax}%
\everymath {\fam0 }
\everydisplay {\fam0 }
  \normalbaselines\bf}

%%%%%%%%%%%%%%%%%%%%%%%%%%%%%%%%%%%%%%%%%%%%%%%%%%%%%%%%%
%% 2) Page set up
\hsize=31pc
\vsize=48pc

\parindent=22pt
\parskip=0pt

\widowpenalty=10000
\clubpenalty=10000

\topskip=12pt 

\skip\footins=20pt
\dimen\footins=3in % maximum footnote height

\abovedisplayskip=6.95pt plus3.5pt minus 3pt
\belowdisplayskip=\abovedisplayskip

%% Output routine

\voffset=7pt\hoffset= .7in%7pt magstep1

\newif\iftitle%!

\def\amheadline{\iftitle%
\hbox to\hsize{\hss\currannalsline\hss}\else\line{\ifodd\pageno
\hfill\thetitle\hfill\llap{\elevenrm\folio}\else\rlap{\elevenrm\folio}
\hfill\theauthors\hfill\fi}\fi}

\headline={\amheadline}%!!!
\footline={\global\titlefalse}
%\output={\bindingoffset\plainoutput}

%%%%%%%%%%%%%%%%%%%%%%%%%%%%%%%%%%%%%%%
% 3) Title page 

 %#1= Volume number, #2=year of publication
\def\annalsline#1#2{\vfill\eject
\ifodd\pageno\else % first page of article on right.
\line{\hfill}
\vfill\eject\fi
\global\titletrue
\def\currannalsline{\eightrm %Annals of Mathematics,%ANNALS
{\eightbf#1} (#2), \thepages}}

\def\titleheadline#1{\def\one{#1}\ifx\one\empty\else
\def\thetitle{{%\frenchspacing%
\let\\ \relax\eightsc\uppercase{#1}}}\fi}

\newif\ifshort

\let\shorttitle\titleheadline

\def\onpages#1#2{\def\thepages{#1--#2}}

\def\thismuchskip[#1]{\vskip#1pt}
\def\ilook{\ifx\next[ \let\go\thismuchskip\else
\let\go\relax\vskip1pt\fi\go}

\def\institution#1{\def\theinstitutions{\vbox{\baselineskip10pt
\def\and{{\eightrm and }}
\def\\{\futurelet\next\ilook}\eightsc #1}}}
\let\institutions\institution

\newwrite\auxfile

\def\startingpage#1{\def\one{#1}\ifx\one\empty\global\pageno=1\else
\global\pageno=#1\fi
\theoremcount=0 \eqcount=0 \sectioncount=0 
\openin1 \jobname.aux \ifeof1 
\onpages{#1}{???}
\else\closein1 \relax\input \jobname.aux
\onpages{#1}{\lastpage}
\fi\immediate\openout\auxfile=\jobname.aux
}

\def\endarticle{\ifRefsUsed\global\RefsUsedfalse%
\else\vskip21pt\theinstitutions%
\nobreak\vskip8pt
%\vbox{\thereceived\therevised}%
\fi%
\write\auxfile{\string\def\string\lastpage{\the\pageno}}}

\outer\def\bye{\endarticle\par \vfill \supereject \end}

% variation on code from amsspt.sty ==>
\def\document{\let\fontlist@\relax\let\alloclist@\relax
 \elevenpoint}%%% add for annals!!!

% <=== end of code varied from amsppt.sty

\newif\ifacks
\long\def\acknowledgements#1{\def\one{#1}\ifx\one\empty\else
\vskip-\baselineskip%
\global\ackstrue\footnote{\ \unskip}{*#1}\fi}

\def\title#1{\titleheadline{#1}
\vbox to80pt{\vfill
\baselineskip=18pt
\parindent=0pt
\overfullrule=0pt
\hyphenpenalty=10000
\everypar={\hskip\parfillskip\relax}
\hbadness=10000
\def\\ {\vskip1sp}
\eighteenbold#1\vskip1sp}}

\newif\ifauthor

\def\author#1{\vskip11pt
\hbox to\hsize{\hss\tenrm By \tensc#1\ifacks\global\acksfalse*\fi\hss}
\ifshort\else\xdef\theauthors{{\eightsc\uppercase{#1}}}\fi%
\vskip21pt\global\authortrue\everypar={\global\authorfalse\everypar={}}}

\def\twoauthors#1#2{\vskip11pt
\hbox to\hsize{\hss%
\tenrm By \tensc#1 {\tenrm and} #2\ifacks\global\acksfalse*\fi\hss}
\ifshort\else\xdef\theauthors{{\eightsc\uppercase{#1 and #2}}}\fi%
\vskip21pt
\global\authortrue\everypar={\global\authorfalse\everypar={}}}

%%%%%%%%%%%%%%%%%%%%%%%%%%%%%%%%
%% 4) Section heads, counters

\newcount\theoremcount
\newcount\sectioncount
\newcount\eqcount

\newif\ifspecialnumon

\def\eqnumber=#1 {\global\eqcount=#1 \global\advance\eqcount by-1\relax}
\def\sectionnumber=#1 {\global\sectioncount=#1 
\global\advance\sectioncount by-1\relax}
\def\proclaimnumber=#1 {\global\theoremcount=#1 
\global\advance\theoremcount by-1\relax}

\newif\ifsection
\newif\ifsubsection

\def\elevenboldmath#1{$#1$\egroup}
\def\mathbold{\hbox\bgroup\elevenbold\elevenboldmath}

\def\section#1{\global\theoremcount=0
\global\eqcount=0
\ifauthor\global\authorfalse\else%
\vskip18pt plus 18pt minus 6pt\fi%
{\parindent=0pt
\everypar={\hskip\parfillskip}%            !!! remove
\def\\ {\vskip1sp}\elevenpoint\bf%
\ifspecialnumon\global\specialnumonfalse$\rm\spnum$%
\gdef\sectnum{$\rm\spnum$}%
\else\interlinepenalty=10000%
\global\advance\sectioncount by1\relax\the\sectioncount%
\gdef\sectnum{\the\sectioncount}%
\fi. \hskip6pt#1%                          !!!add }} and stop here
\vrule width0pt depth12pt}
\hskip\parfillskip%\break%!
\global\sectiontrue%
\everypar={\global\sectionfalse\global\interlinepenalty=0\everypar={}}%
\ignorespaces

}

%%%%%%%%%%%%%%%%%%%%%%%%%%%%%%%%
%% 5) Equation Macros

\newif\ifspequation

\let\eqno\leqno %automatic left side equation numbers %%!!!remove l-eqno

\newif\ifineqalignno
\let\saveleqalignno\leqalignno                        %%!!!remove l-eqno
\def\leqalignno{\let\eqnu\Eeqnu\saveleqalignno}

\let\eqalignno\leqalignno

\def\sectandeqnum{%
\ifspecialnumon\global\specialnumonfalse
$\rm\spnum$\gdef\eqnum{{$\rm\spnum$}}\else\global\firstlettertrue
\global\advance\eqcount by1 
\ifappend\applett\else\the\sectioncount\fi.%
\the\eqcount
\xdef\eqnum{\ifappend\applett\else\the\sectioncount\fi.\the\eqcount}\fi}

\def\eqnu{\leqno{\hbox{\elevenrm\ifspequation\else(\fi\sectandeqnum
\ifspequation\global\spequationfalse\else)\fi}}}      %!!! l-eqno

\def\Speqnu{\global\setbox\leqnobox=\hbox{\elevenrm
\ifspequation\else%
(\fi\sectandeqnum\ifspequation\global\spequationfalse\else)\fi}}

\def\Eeqnu{\hbox{\elevenrm
\ifspequation\else%
(\fi\sectandeqnum\ifspequation\global\spequationfalse\else)\fi}}

\newif\iffirstletter
\global\firstlettertrue
\def\eqletter#1{\global\specialnumontrue\iffirstletter\global\firstletterfalse
\global\advance\eqcount by1\fi
\gdef\spnum{\the\sectioncount.\the\eqcount#1}\eqnu}

%%% Split math
\newbox\leqnobox
\def\outsideeqnu#1{\global\setbox\leqnobox=\hbox{#1}}

\def\eatone#1{}

%% Vertically centers equation number.
\def\dosplit#1#2{\vskip-.5\abovedisplayskip
\setbox0=\hbox{$\let\eqno\outsideeqnu%
\let\eqnu\Speqnu\let\leqno\outsideeqnu#2$}%
\setbox1\vbox{\noindent\hskip\wd\leqnobox\ifdim\wd\leqnobox>0pt\hskip1em\fi%
$\displaystyle#1\mathstrut$\hskip0pt plus1fill\relax
\vskip1pt
\line{\hfill$\let\eqnu\eatone\let\leqno\eatone%
\displaystyle#2\mathstrut$\ifmathqed~~\qed\fi}}%
\copy1
\ifvoid\leqnobox
\else\dimen0=\ht1 \advance\dimen0 by\dp1
\vskip-\dimen0
\vbox to\dimen0{\vfill
\hbox{\unhbox\leqnobox}
\vfill}
\fi}

\everydisplay{\lookforbreak}

\long\def\lookforbreak #1$${\def\mathone{#1}
\expandafter\testforbreak\mathone\splitmath @}

\def\testforbreak#1\splitmath #2@{\def\mathtwo{#2}\ifx\mathtwo\empty%
#1$$%
\ifmathqed\vskip-\belowdisplayskip
\setbox0=\vbox{\let\eqno\relax\let\eqnu\relax$\displaystyle#1$}% 
\vskip-\ht0\vskip-3.5pt\hbox to\hsize{\hfill\qed}
\vskip\ht0\vskip3.5pt\fi
\else$$\vskip-\belowdisplayskip
\vbox{\dosplit{#1}{\let\eqno\eatone
\let\splitmath\relax#2}}%
\nobreak\vskip.5\belowdisplayskip
\noindent\ignorespaces\fi}

%% Proof box to be used when proof ends with equation.

\newif\ifmathqed

%%%%%%%%%%%%%%%%%%%%%%%%%%%%%
%% \mtable, Math table to make binary table easily

%% Use:
% \mtable
% &n_1&n_2&n_3&n_4&n_5&n_6\cr
% \Delta_1&M_3&M_2&0&0&0&0\cr
% \Delta_2&0&0&M_1&M_3&0&0\cr
% \endmtable

\newcount\linenum
\newcount\colnum

%++
\def\spline{\omit&\multispan{\the\colnum}{\hrulefill}\cr}
\def\colcounter{\ifnum\linenum=1\global\advance\colnum by1\fi}

\def\everyline{\noalign{\global\advance\linenum by1\relax}%
\ifnum\linenum=2\spline\fi}

\def\mtable{\bgroup\offinterlineskip
\everycr={\everyline}\global\linenum=0
\halign\bgroup\vrule height 10pt depth 4pt width0pt
\hfill$##$\hfill\hskip6pt\ifnum\linenum>1
\vrule\fi&&\colcounter\hskip12pt\hfill$##$\hfill\hskip12pt\cr}

\def\endmtable{\crcr\egroup\egroup}

%%%%%%%%%%%%%%%%%%%%%%%%%%%%%
% Array

%% Will work in math or in text, will be in math mode inside array. 
%% For each column desired supply
%% r, l, or c, for right, left, or center orientation of that column.
%% End each line with \\.

%% To use: 
%  \array ccc*
%  x_s\leq a_1\\
%  a_s<x_s^s<b_s\\
%  x_s\geq a_1
%  \endarray

\def\xast{*}
\newcount\intable
\newcount\mathcol
\newcount\savemathcol
\newcount\topmathcol
\newdimen\arrayhspace
\newdimen\arrayvspace

\arrayhspace=8pt % horizontal space between columns, (half this width
                 %  will horizontally precede and follow the array)
\arrayvspace=12pt % vertical space between lines

\newif\ifdollaron

\def\mathalign#1{\def\arg{#1}\ifx\arg\xast%
\let\go\relax\else\let\go\mathalign%
\global\advance\mathcol by1 %
\global\advance\topmathcol by1 %
\expandafter\def\csname  mathcol\the\mathcol\endcsname{#1}%
\fi\go}

\def\arraypickapart#1]#2*{\if#1c \ifmmode\vcenter\else
\global\dollarontrue$\vcenter\fi\else%
\if#1t\vtop\else\if#1b\vbox\fi\fi\fi\bgroup%
\def\one{#2}}

\def\arraystrut{\vrule height .7\arrayvspace depth .3\arrayvspace width 0pt}

\def\array#1#2*{\def\firstarg{#1}%
\if\firstarg[ \def\two{#2} \expandafter\arraypickapart\two*\else%
\ifmmode\vcenter\else\vbox\fi\bgroup \def\one{#1#2}\fi%
\global\everycr={\noalign{\global\mathcol=\savemathcol\relax}}%
\def\\ {\cr}%
\global\advance\intable by1 %
\ifnum\intable=1 \global\mathcol=0 \savemathcol=0 %
\else \global\advance\mathcol by1 \savemathcol=\mathcol\fi%
\expandafter\mathalign\one*%
\mathcol=\savemathcol %
\halign\bgroup&\hskip.5\arrayhspace\arraystrut%
\global\advance\mathcol by1 \relax%
\expandafter\if\csname mathcol\the\mathcol\endcsname r\hfill\else%
\expandafter\if\csname mathcol\the\mathcol\endcsname c\hfill\fi\fi%
$\displaystyle##$%
\expandafter\if\csname mathcol\the\mathcol\endcsname r\else\hfill\fi\relax%
\hskip.5\arrayhspace\cr}

\def\endarray{\crcr\egroup\egroup%
\global\mathcol=\savemathcol %
\global\advance\intable by -1\relax%
\ifnum\intable=0 %
\ifdollaron\global\dollaronfalse $\fi
\loop\ifnum\topmathcol>0 %
\expandafter\def\csname  mathcol\the\topmathcol\endcsname{}%
\global\advance\topmathcol by-1 \repeat%
\global\everycr={}\fi%
}

\def\big#1{{\hbox{$\left#1\vbox to 10pt{}\right.\n@space$}}}
\def\Big#1{{\hbox{$\left#1\vbox to 13pt{}\right.\n@space$}}}
\def\bigg#1{{\hbox{$\left#1\vbox to 16pt{}\right.\n@space$}}}
\def\Bigg#1{{\hbox{$\left#1\vbox to 19pt{}\right.\n@space$}}}

%%%%%%%%%%%%%%%%%%%%%%%%%%%%%%%%%%%%%%%%%%%%%%%%%%%%%%%%%%%%%%%%
% 6) Figure and Table Captions.

\def\figcaption#1#2#3{\topinsert
\vskip4pt %<===topadjust to match height of ascenders on opposing page.
\vbox to#3{\vfill}\vskip1sp
\setbox0=\hbox{\eightsc Figure #1.\hskip12pt\eightpoint #2}
\ifdim\wd0>\hsize
\noindent\eightsc Figure #1.\hskip12pt\eightpoint #2
\else
\centerline{\eightsc Figure #1.\hskip12pt\eightpoint #2}
\fi
\vskip16pt
\endinsert}

\def\wfig#1#2#3{\topinsert
\vskip4pt %<===topadjust to match height of ascenders on opposing page.
\hbox to\hsize{\hss\vbox{\hrule height .25pt width #3
\hbox to #3{\vrule width .25pt height #2\hfill\vrule width .25pt height #2}
\hrule height.25pt}\hss}
\vskip1sp
\centerline{\eightsc Figure #1}
\vskip16pt
\endinsert}

\def\wfigcaption#1#2#3#4{\topinsert
\vskip4pt %<===topadjust to match height of ascenders on opposing page.
\hbox to\hsize{\hss\vbox{\hrule height .25pt width #4
\hbox to #4{\vrule width .25pt height #3\hfill\vrule width .25pt height #3}
\hrule height.25pt}\hss}
\vskip1sp
\setbox0=\hbox{\eightsc Figure #1.\hskip12pt\eightpoint\rm #2}
\ifdim\wd0>\hsize
\noindent\eightsc Figure #1.\hskip12pt\eightpoint\rm #2\else
\centerline{\eightsc Figure #1.\hskip12pt\eightpoint\rm #2}\fi
\vskip16pt
\endinsert}

\def\tabcaption#1#2{\vskip6pt
\setbox0=\hbox{\eightsc Table #1.\hskip12pt\eightpoint #2}
\ifdim\wd0>\hsize
\noindent\eightsc Table #1.\hskip12pt\eightpoint #2
\else
\centerline{\eightsc Table #1.\hskip12pt\eightpoint #2}
\fi
\vskip6pt}

\def\endinsert{\egroup\if@mid\dimen@\ht\z@\advance\dimen@\dp\z@ 
\advance\dimen@ 12\p@\advance\dimen@\pagetotal\ifdim\dimen@ >\pagegoal 
\@midfalse\p@gefalse\fi\fi\if@mid\smallskip\box\z@\bigbreak\else
\insert\topins{\penalty 100 \splittopskip\z@skip\splitmaxdepth\maxdimen
\floatingpenalty\z@\ifp@ge\dimen@\dp\z@\vbox to\vsize {\unvbox \z@ 
\kern -\dimen@ }\else\box\z@\nobreak\smallskip\fi}\fi\endgroup}

\def\pagecontents{
\ifvoid\topins \else\iftitle\else 
\unvbox \topins \fi\fi \dimen@ =\dp \@cclv \unvbox 
\@cclv 
\ifvoid\topins\else\iftitle\unvbox\topins\fi\fi
\ifvoid \footins \else \vskip \skip \footins \footnoterule 
\unvbox \footins \fi \ifr@ggedbottom \kern -\dimen@ \vfil \fi}

%%%%%%%%%%%%%%%%%%%%%%%%%%%%%%%%%%%%%%%%%%%%%%%%%%%%%%%%%%%%%%%%
% 7) End Matter

\newif\ifappend

\def\appendix#1#2{\def\applett{#1}\def\two{#2}%
\global\appendtrue
\global\theoremcount=0
\global\eqcount=0
\vskip18pt plus 18pt
\vbox{\parindent=0pt
\everypar={\hskip\parfillskip}
\def\\ {\vskip1sp}\elevenbold Appendix% 
\ifx\applett\empty\gdef\applett{A}\ifx\two\empty\else.\fi%
\else\ #1.\fi\hskip6pt#2\vskip12pt}%
\global\sectiontrue%
\everypar={\global\sectionfalse\everypar={}}\nobreak\ignorespaces}

\newif\ifRefsUsed
\long\def\references{\global\RefsUsedtrue\vskip21pt%
\theinstitutions
\global\everypar={}\global\bibnum=0
\vskip20pt\goodbreak\bgroup%
\vbox{\centerline{\eightsc References}\vskip6pt}% 
\ifdim\maxbibwidth>0pt
\leftskip=\maxbibwidth%
\parindent=-\maxbibwidth%
\else
\leftskip=18pt%
\parindent=-18pt%
\fi
\ninepoint%
\frenchspacing
\nobreak\ignorespaces\everypar={\amref}%
}

\def\endreferences{\vskip1sp\egroup\global\everypar={}%
\nobreak\vskip8pt\vbox{\thereceived\therevised}
}

\newcount\bibnum

\def\amref#1 {\global\advance\bibnum by1%
\immediate\write\auxfile{\string\expandafter\string\def\string\csname
\space #1croref\string\endcsname{[\the\bibnum]}}%
\leavevmode\hbox to18pt{\hbox to13.2pt{\hss[\the\bibnum]}\hfill}}

\def\bibline{\hbox to30pt{\hrulefill}\/\/}

\def\name#1{{\eightsc#1}}

\newdimen\maxbibwidth
\def\AuthorRefNames [#1] {%
\immediate\write\auxfile{\string\def\string\cite\string##1{[\string##1]}}

\def\amref{\spamref}
\setbox0=\hbox{[#1] }\global\maxbibwidth=\wd0\relax}

\def\spamref[#1] {\leavevmode\hbox to\maxbibwidth{\hss[#1]\hfill}}

%%%%%%%%%%%%%%%%%%%%%%%%%%%%%%%%%%%%%%%%%%%%%%%%%%%%%%%%%%%%%%%%
%% 8) Footnotes

\def\footnoterule{\kern-3pt\hrule width1in height.5pt\kern2.5pt}

\def\footnote#1#2{%
\plainfootnote{#1}{{\eightpoint\normalbaselineskip11pt
\normalbaselines#2}}}

\def\vfootnote#1{%
\insert \footins \bgroup \eightpoint\baselineskip11pt
\interlinepenalty \interfootnotelinepenalty
\splittopskip \ht \strutbox \splitmaxdepth \dp \strutbox \floatingpenalty 
\@MM \leftskip \z@skip \rightskip \z@skip \spaceskip \z@skip 
\xspaceskip \z@skip
{#1}$\,$\footstrut \futurelet \next \fo@t}

%%%%%%%%%%%%%%%%%%%%%%%%%%%%%%%%%%%%%%%%%%%%%%%%%%%%%%%%%%%%%%%%
%% 9) Theorem type environments

\newif\iffirstadded
\newif\ifadded

\def\addedlett{}

\def\alltheoremnums{%
\ifspecialnumon\global\specialnumonfalse
\ifadded\global\addedfalse
\iffirstadded\global\firstaddedfalse
\global\advance\theoremcount by1 \fi
\ifappend\applett\else\the\sectioncount\fi.\the\theoremcount\addedlett%
\xdef\theoremnum{\ifappend\applett\else\the\sectioncount\fi.%
\the\theoremcount\addedlett}%
\else$\rm\spnum$\def\theoremnum{{$\rm\spnum$}}\fi%
\else\global\firstaddedtrue
\global\advance\theoremcount by1 
\ifappend\applett\else\the\sectioncount\fi.\the\theoremcount%
\xdef\theoremnum{\ifappend\applett\else\the\sectioncount\fi.%
\the\theoremcount}\fi}

\def\allcorolnums{%
\ifspecialnumon\global\specialnumonfalse
\ifadded\global\addedfalse
\iffirstadded\global\firstaddedfalse
\global\advance\corolcount by1 \fi
\the\corolcount\addedlett%
\else$\rm\spnum$\def\corolnum{$\rm\spnum$}\fi%
\else\global\advance\corolcount by1 
\the\corolcount\fi}

%% use for Theorem, Corollary, Lemma, Proposition, Demonstration and similar.

%\newcount\corolcount
%\def\xcorol{Corollary} 
%\def\xtheorem{Theorem}
%\def\xmaintheorem{Main Theorem} %!!! -was changed to numerate

\newif\ifthtitle

\let\saverparen)
\let\savelparen(
\def\rmparenl{{\rm(}}
\def\rmparenr{{\rm\/)}}
{
\catcode`(=13
\catcode`)=13
\gdef\makeparensRM{\catcode`(=13\catcode`)=13\let(=\rmparenl%
\let)=\rmparenr%
\everymath{\let(\savelparen%
\let)\saverparen}%
\everydisplay{\let(\savelparen%
\let)\saverparen\lookforbreak}}}

\medskipamount=8pt plus.1\baselineskip minus.05\baselineskip

\def\rmtext#1{\hbox{\rm#1}}

\def\proclaim#1{\vskip-\lastskip
\def\one{#1}\ifx\one\xtheorem\global\corolcount=0\fi
\ifsection\global\sectionfalse\vskip-6pt\fi
\medskip
{\elevensc#1}%
\ifx\one\xmaintheorem\global\corolcount=0
\gdef\theoremnum{Main Theorem}\else%
\ifx\one\xcorol\ \allcorolnums\else\ \alltheoremnums\fi\fi%
\ifthtitle\ \global\thtitlefalse{\rm(\thethtitle)}\fi.%
\hskip1em\bgroup\let\text\rmtext\makeparensRM\it\ignorespaces}

\def\nonumproclaim#1{\vskip-\lastskip
\def\one{#1}\ifx\one\xtheorem\global\corolcount=0\fi
\ifsection\global\sectionfalse\vskip-6pt\fi
\medskip
{\elevensc#1}.\ifx\one\xmaintheorem\global\corolcount=0
\gdef\theoremnum{Main Theorem}\fi\hskip.5pc%
\bgroup\it\makeparensRM\ignorespaces}

\def\endproclaim{\egroup\medskip}

%% Use demo for Proof, Proof of, Definition, Example,
%% Remark, Case, Subcase, Conjecture, Note, Notation,
%% Convention, Construction and Step.
%% Any other use for demo will format similar to `Proof.'

\def\xproof{Proof}
\def\xremark{Remark}
\def\xcase{Case}
\def\xsubcase{Subcase}
\def\xconjecture{Conjecture}
\def\xstep{Step}
\def\xof{of}

\def\deconstruct#1 #2 #3 #4 #5 @{\def\one{#1}\def\two{#2}\def\three{#3}%
\def\four{#4}%
\ifx\two\empty #1\else%
\ifx\one\xproof%
\ifx\two\xof%
  \ifx\three\xcorol Proof of Corollary \rm#4\else%
     \ifx\three\xtheorem Proof of Theorem \rm#4\else\xone\fi%
  \fi\fi%
\else\xone\fi\fi.}

\def\pickup#1 {\def\this{#1}%
\ifx\this\xproof\global\let\go\demoproof
\global\let\enddemo\endproof\else
\ifx\this\xremark\global\let\go\demoremark\else
\ifx\this\xcase\global\let\go\demostep\else
\ifx\this\xsubcase\global\let\go\demostep\else
\ifx\this\xconjecture\global\let\go\demostep\else
\ifx\this\xstep\global\let\go\demostep\else
\global\let\go\demoproof\fi\fi\fi\fi\fi\fi}

\newif\ifnonum
\def\demo#1{\vskip-\lastskip
\ifsection\global\sectionfalse\vskip-6pt\fi
\def\one{#1 }\def\two{#1*}%
\setbox0=\hbox{\expandafter\pickup\one}\expandafter\go\two}

\def\numbereddemo#1{\vskip-\lastskip
\ifsection\global\sectionfalse\vskip-6pt\fi
\def\two{#1*}%
\expandafter\demoremark\two}

\def\demoproof#1*{\medskip\def\xone{#1}
{\ignorespaces\it\expandafter\deconstruct\xone {} {} {} {} {} @%
\unskip\hskip6pt}\rm\ignorespaces}

\def\demoremark#1*{\medskip
{\it\ignorespaces#1\/} \ifnonum\global\nonumtrue\else
 \alltheoremnums\unskip.\fi\hskip1pc\rm\ignorespaces}

\def\demostep#1 #2*{\vskip4pt
{\it\ignorespaces#1\/} #2.\hskip1pc\rm\ignorespaces}

\def\enddemo{\medskip}

\def\endproof{\ifmathqed\global\mathqedfalse\medskip\else
\parfillskip=0pt~~\hfill\qed\medskip
\fi\global\parfillskip0pt plus 1fil\relax
\gdef\enddemo{\medskip}}

\def\qed{\vbox{\hrule\hbox{\vrule height6pt\hskip6pt\vrule}\hrule}}

%% Proof box to be used in a \proclaim{}...\endproclaim environment

\def\proofbox{\parfillskip=0pt~~\hfill\qed\vskip1sp\parfillskip=
0pt plus 1fil\relax}

%%%%

%%%%%%%%%%%%%%%%%%%%%
%% 10) CrossRefs 

%%% Generic crossreferencing
%%% to use: \label\nameoflabel* (will give the page number when referenced)

% Commands to access current state of counter, for cross-referencing
% \sectnum
% \theoremnum
% \eqnum

%%% You can make another definition that includes counters and/or the
%%% page number and access this information as the second argument:
%%% \label\yourlabelname[2.13]*

%%% Since this method of cross-referencing relies
%%% on an auxiliary file, the first time you tex the file
%%% you will get `??' when you write \ref\nameoflabel.
%%% When you TeX the file the second time the auxiliary file
%%% will be input and \ref\nameoflabel will produce the cross-ref.

\def\stripbs#1#2*{\def\one{#2}}

\def\emptyspace{ }
\def\nextthing{}
\def\newline{***}
\def\eatone#1{ }

\def\lookatspace#1{\ifcat\noexpand#1\ \else%
\gdef\nextthing{}\xdef\next{#1}%
\ifx\next\emptyspace%
\let\nextthing\emptyspace\else\ifx\next\newline%
\gdef\nextthing{\eatone}\fi\fi\fi\egroup\nextthing#1}

{\catcode`\^^M=\active%
\gdef\spacer{\bgroup\catcode`\^^M=\active%
\let^^M=\newline\obeyspaces\lookatspace}}

\def\ref#1{\seeifdefined{#1}\expandafter\csname\one\endcsname\spacer}

\def\cite#1{\expandafter\ifx\csname#1croref\endcsname\relax[??]\else
\csname#1croref\endcsname\fi\spacer}

%% for testing in \label and \ref to see if term already labeled.

\def\seeifdefined#1{\expandafter\stripbs\string#1croref*%
\crorefdefining{#1}}

\newif\ifcromessage
\global\cromessagetrue

\def\crorefdefining#1{\ifdefined{\one}{}
{\ifcromessage\global\cromessagefalse%
\message{\spaces\spaces\spaces\spaces\spaces\spaces\spaces}%
\message{<Undefined reference.}%
\message{Please TeX file once more to have accurate cross-references.>}%
\message{\spaces\spaces\spaces\spaces\spaces\spaces\spaces}\fi[??]}}

\def\label#1#2*{\gdef\ctest{#2}%
\xdef\currlabel{\string#1croref}
\expandafter\seeifdefined{#1}%
\ifx\empty\ctest%
\xdef\labelnow{\write\auxfile{\noexpand\def\currlabel{\the\pageno}}}%
\else\xdef\labelnow{\write\auxfile{\noexpand\def\currlabel{#2}}}\fi%
\labelnow}

\def\ifdefined#1#2#3{\expandafter\ifx\csname#1\endcsname\relax%
#3\else#2\fi}

%%%%%%%%%%%%%%%%%%%%%
%% 11) Listing

%% To use with asterisks:

%%%%%%%%%%%%%%%%%%%%%%
%% 12) Article and Journal Table of Contents

\def\articlecontents{
\vskip20pt\centerline{\bf Table of Contents}\everypar={}\vskip6pt
\bgroup \leftskip=3pc \parindent=-2pc 
\def\item##1{\vskip1sp\indent\hbox to2pc{##1.\hfill}}}

\def\endcontents{\vskip1sp\leftskip=0pt\egroup}

\def\journalcontents{\vfill\eject
\def\currannalsline{\hfill}
\global\titletrue
\vglue3.5pc
\centerline{\tensc\hskip12pt TABLE OF CONTENTS}\everypar={}\vskip30pt
\bgroup \leftskip=34pt \rightskip=-12pt \parindent=-22pt 
  \def\\ {\vskip1sp\noindent}
\def\pagenum##1{\unskip\parfillskip=0pt\dotfill##1\vskip1sp
\parfillskip=0pt plus 1fil\relax}
\def\name##1{{\tensc##1}}}

%% default values

\institution{}
\onpages{0}{0}
\def\lastpage{???}
\def\thetitle{Title ???}
\def\theauthors{Authors ???}
\def\thereceived{}
\def\therevised{}

\gdef\split{\relaxnext@\ifinany@\let\next\insplit@\else
 \ifmmode\ifinner\def\next{\onlydmatherr@\split}\else
 \let\next\outsplit@\fi\else
 \def\next{\onlydmatherr@\split}\fi\fi\let\eqnu\xspliteqnu\next}

\gdef\align{\relaxnext@\ifingather@\let\next\galign@\else
 \ifmmode\ifinner\def\next{\onlydmatherr@\align}\else
 \let\next\align@\fi\else
 \def\next{\onlydmatherr@\align}\fi\fi\let\eqnu\xspliteqnu\next}

\def\spliteqnu{{\tenrm\sectandeqnum}\relax}

\def\xspliteqnu{\tag\spliteqnu}

\catcode`@=12

\document

%-------------- Publisher's entries --------------------

\annalsline{June}{1996}
%\line{\hfil Revised }
\line{\ninerm \hfill Dedicated to Yuri I. Manin}
\vskip -0.2cm
\line {\ninerm \hfill on the occasion of his 60 birthday}
\startingpage{1}     %numeration
%%\received{??}
%\revised{}

%\magnification=\magstep1

%--------------- Author macros ---------------
%                   MACROS
%
%                                 AUX
%
%
%
%                      endaux
%

\def\iif{\quad\hbox{ if }\quad}

\def\for{\  \hbox{ for } \ }
\def\if{ \ \hbox{ if } \ }

\def\where{\  \hbox{ where } \ }
\def\and{\  \hbox{ and } \ }

\def\equal{\buildrel  def \over =}
\def\lan{\langle}
\def\ran{\rangle}
\def\llan{\langle\!\langle}
\def\rran{\rangle\!\rangle}

\def\la{\lambda}

\def\om{\omega}

\def\th{\theta}
\def\al{\alpha}
\def\be{\beta}
\def\ga{\gamma}
\def\ep{\epsilon}

\def\de{\delta}
\def\De{\Delta}
\def\ka{\kappa}
\def\si{\sigma}
\def\Si{\Sigma}
\def\Ga{\Gamma}
\def\ze{\zeta}

    %from copy, ell
\def\pa{\partial}

\def\vph{\varphi}

\def\vep{\varepsilon}

\def\tal{\tilde{\alpha}}

\def\txi{\tilde{\xi}}

\def\tU{\tilde{U}}
\def\tw{\tilde w}
\def\tW{\tilde W}

\def\tz{\tilde z}
\def\tb{\tilde b}
\def\ta{\tilde a}

\def\ty{\tilde y}

\def\tG{\tilde G}

\def\tPhi{\tilde {\Phi}}

\def\hH{\hat{H}}

\def\hT{\hat{T}}

\def\hw{\hat{w}}
\def\hu{\hat{u}}

\def\hv{\hat{v}}

\def\he{\hat{e}}
\def\hf{\hat{f}}
\def\hde{\hat{\delta}}

\def\C{\bold{C}}
\def\Q{\bold{Q}}

\def\R{\bold{R}}
\def\N{\bold{N}}
\def\Z{\bold{Z}}

\def\one{\bold{1}}

\def\0{\bold{0}}

\def\C{\hbox{\bf C}}

\def\H{\bold{H}}% macdonald

\def\f{\Cal{F}}

\def\r{\Cal{R}}
\def\l{\Cal{L}}

\def\k{\Cal{K}}

\def\p{\Cal{P}}
\def\a{\Cal{A}}
\def\h{\Cal{H}}

\def\v{\Cal{V}}

\def\o{\Cal{O}}

\def\Indcx{\mathop{\hbox{\rm Ind}\,}_{\hbox{\eightbf C}[X]}\,}

\font\germ=eufb10 %at 12pt 
%\font\germm=eufb9 at 12pt                       
%\font\germ=eufm9 at 12pt 
\def\goth#1{\hbox{\germ #1}}

\def\TT{\goth{T}}
\def\HH{\goth{H}}

\def\AA{\goth{A}}

\font\smm=msbm10 at 12pt 
\def\symbol#1{\hbox{\smm #1}}
\def\lsmash{{\symbol n}}

%endmacros

%------------------------------------------------------------------
%-------------- Author entries --------------------

%\comment                               %to remove the title
\title
{Intertwining operators \\
of double affine Hecke algebras}
 
 %Article title
\shorttitle{ Intertwiners of double Hecke algebras}
 % Shortened version for headline title

% Acknowledgements: Please enter all acknowledgements here.
\acknowledgements{
Partially supported by NSF grant DMS--9622829 and NWO, the Netherlands}

% Please uncomment and use appropriate command:
\author{ Ivan Cherednik}
%\twoauthors{}{}
%\authors{}% Separate each author with a comma and a space.

% Institution:
\institutions{
Math. Dept, University of North Carolina at Chapel Hill,   
 N.C. 27599-3250
\\ Internet: chered\@math.unc.edu
}

%\endcomment                              %to remove the title
%-------------- Article Text--------------------
%\intro %(Optional, Introduction)
%
%
%			
%                        INTRO
%
%
%{\bf 0. Introduction.
\vfil%%

Continuing [C3,C4], we study the intertwining
operators of double affine Hecke algebras $\HH$. They appeared in several
papers (especially in [C2,C4,C6]). However  for the first time here
we apply them systematically to create the nonsymmetric  [M3,C4] 
and symmetric [M2] Macdonald  polynomials 
for arbitrary root systems and to start the theory of induced and
co-spherical $\HH$-modules.

The importance of this technique was clearly demonstrated
in  recent papers by F. Knop and S. Sahi [Kn],[KS],[S]. Using the
intertwiners of the double affine Hecke algebras 
in the case of $GL$ (dual to those considered in [C1,C2]) 
they proved the $q,t$-integrality
conjecture by I. Macdonald [M1] and managed to establish the
positivity of the coefficients of the Macdonald polynomials
in the differential case. 
As to the integrality, we mention another approach based on the 
so-called Vinet operators (see [LV] and a recent work by Kirillov,
Noumi), and the results by Garsia, Remmel, and Tesler.

We do not try in this paper
to get the best possible estimates for the denominators of the
Macdonald polynomials (generally speaking, the problem
looks more complicated than in the stable $GL$-case). However
even rather straightforward analysis of the intertwiners
gives a lot. For instance, it is enough 
to ensure the existence of the restricted Macdonald polynomials  
at roots of unity from [C3,C4], where we
used less convenient methods based directly on the definition
or on the recurrence relations. 

The technique of intertwiners combined with the (projective)
action of
$GL(2,\Z)$ from [C3] gives another proof
of the norm and the evaluation formulas (see [C4]).
%To be more exact, mostly we use the automorphism of $\HH$ corresponding 
%to $\Bigl(\matrix -1&  0\\  1& 1 \endmatrix \Bigr)$.
Here  the $\HH$-embedding of the space of
nonsymmetric polynomials into the space of
functions on the affine Weyl
group $\tW$ ([C4], Proposition 5.2) plays a key role. 
The latter  representation  when restricted to the affine Hecke 
subalgebra turns into the classical one from [IM]
as $t$ is a power of $p$ and $q\to 0$
($\tW$ is identified with the set of double cosets of the corresponding 
$p$-adic group with respect to the Iwahori subroup).

Another important application is a calculation of
the Fourier transforms 
of the Macdonald polynomials
in the sense of [C3,C4].  For instance, it gives a
canonical  identification of
the polynomial representation of the affine Hecke algebra  
with the representation in functions on the weight lattice
(which collapses in the $p$-adic limit).

We introduce  a proper discretization
of the  $\mu$-function (the truncated theta-function
making Macdonald's polynomials
pairwise orthogonal) and the corresponding discrete
inner product on $\hbox{Funct}(\tW)$.
It readily gives the proportionality of
the norms of the Macdonald polynomials [M2,C2,M3,C4]
and those defined
for the  Jackson integral taken instead of the
constant term in the inner product. The 
coefficient of proportionality is described by the Aomoto
conjecture (see  [A,Ito]) recently proved by Macdonald 
(to calculate it one can also follow  [C2],
 replacing the shift operators by  their discretizations).

We note that the Macdonald polynomials considered as functions 
on $\tW$ are square integrable for finitely many
weights only. Here  $|q|\neq 1$ and 
the real part $\Re(k)$ for $t=q^k$
is to be negative (otherwise we  have none).
The program is to describe all 
integrable and non-integrable 
eigenfunctions of the discrete Dunkl operators
in this representation and to study
the corresponding Fourier transform.
In contrast
to the classical $p$-adic harmonic analysis (see e.g. [HO]) 
the Plancherel measure coincides with the discretization of $\mu$ 
(the Fourier transform is self-dual). 

More generally, 
we  consider the action of the double affine Hecke algebra in
the same space $\hbox{Funct}(\tW)$ 
depending on an arbitrary given weight.
Its submodule generated by the delta-functions
is  induced (from a character of the standard 
polynomial subalgebra) and co-spherical.
Mainly following [C5], we find out when  arbitrary induced
representations (in the same sense)  are irreducible and
co-spherical using the technique of intertwiners.
The answer is a natural "affinization" of
the  well-known statements in the
$p$-adic case (see e.g. [KL], [C5]). 
The classification of co-spherical 
representations is important for the harmonic analysis and
plays the key role in the theory
of affine Knizhnik-Zamolodchikov equations 
(see [C6,C7,C8]). We also 
induce up irreducible representations of affine Hecke subalgebras
([C6] is devoted to applications of such representations).
If $q$ is sufficiently general the \HH\ -modules we get
are irreducible,
so one can use the classification of [KL].

Thus in this paper we begin a systematic 
study of the representations of double affine
Hecke algebras and related harmonic analysis.
The polynomial representation
considered in the series of papers [C2-4] 
devoted to the  Macdonald conjectures is remarkable,
but still just an  example.

\vfil%%%
The paper was started during my stay at RIMS
(Kyoto University), continued at CRM in Montreal,
and completed at the University of Nijmegen.
I am grateful to T. Miwa, L. Vinet, G. Heckman
and my colleagues at these institutes
for the kind invitations and the hospitality.
The author thanks E. Frenkel, G. Heckman,
D. Kazhdan, I. Macdonald, and E. Opdam
for useful discussions.

%
%
%		Section 1
%
%
%\vskip 10pt
\section {Affine Weyl groups}

Let $R=\{\al\}   \subset \R^n$ be a root system of type $A,B,...,F,G$
with respect to a euclidean form $(z,z')$ on $\R^n \ni z,z'$,
normalized by the standard condition
that $(\al,\al)=2$ for long $\al$.
Let us  fix the set $R_{+}$ of positive  roots ($R_-=-R_+$), 
the corresponding simple 
roots $\al_1,...,\al_n$, and  their dual counterparts 
$a_1 ,..., a_n,  a_i =\al_i^\vee, \where \al^\vee =2\al/(\al,\al)$.  
The dual fundamental weights
$b_1,...,b_n$  are determined from the relations  $ (b_i,\al_j)= 
\de_i^j $ for the 
Kronecker delta. We will also use the dual root system
$R^\vee =\{\al^\vee, \al\in R\}, R^\vee_+$, and the lattices
$$
\eqalignno{
& A=\oplus^n_{i=1}\Z a_i \subset B=\oplus^n_{i=1}\Z b_i, 
}
$$
  $A_\pm, B_\pm$  for $\Z_{\pm}=\{m\in\Z, \pm m\ge 0\}$
instead of $\Z$. (In the standard notations, $A= Q^\vee,\ 
B = P^\vee $ - see [B].)  Later on,
$$
\eqalign{
&\nu_{\al}=\nu_{\al^\vee}\ =\ (\al,\al),\  \nu_i\ =\ \nu_{\al_i}, \ 
\nu_R\ = \{\nu_{\al}, \al\in R\}\subset \{2,1,2/3\}.
}
$$
$$\eqalign{
&\rho_\nu\ =\ (1/2)\sum_{\nu_{\al}=\nu} \al \ =
\ (\nu/2)\sum_{\nu_i=\nu}  b_i, \for\al\in R_+.
%\cr
%&r_\nu\ =\ \rho_\nu^\vee \ =\ (2/\nu)\rho_\nu\ =\ 
%\sum_{\nu_i=\nu}  b_i,\quad 2/\nu=1,2,3. 
} 
\eqnu
\label\rhor\eqnum*
$$

The vectors $\ \tal=[\al,k] \in 
\R^n\times \R \subset \R^{n+1}$ 
for $\al \in R, k \in \Z $ 
form the {\it affine root system} 
$R^a \supset R$ ( $z\in \R^n$ are identified with $ [z,0]$).
We add  $\al_0 \equal [-\th,1]$ to the  simple roots 
for the {\it maximal root} $\th \in R$.
The corresponding set $R^a_+$ of positive roots coincides
with $R_+\cup \{[\al,k],\  \al\in R, \  k > 0\}$. 

We denote the Dynkin diagram and its affine completion with
$\{\al_j,0 \le j \le n\}$ as the vertices by $\Ga$ and $\Ga^a$
($m_{ij}=2,3,4,6$\  if $\al_i\and\al_j$ are joined by 0,1,2,3 laces
respectively).
The set of
the indices of the images of $\al_0$ by all 
the automorphisms of $\Ga^a$ will be denoted by $O$ ($O=\{0\} 
\for E_8,F_4,G_2$). Let $O^*={r\in O, r\neq 0}$.
The elements $b_r$ for $r\in O^*$ are the so-called minuscule
weights ($(b_r,\al)\le 1$ for
$\al \in R_+$).

Given $\tal=[\al,k]\in R^a,  \ b \in B$, let  
$$
\eqalignno{
&s_{\tal}(\tz)\ =\  \tz-(z,\al^\vee)\tal,\ 
\ b'(\tz)\ =\ [z,\ze-(z,b)]
&\eqnu
\label\saction\eqnum*
}
$$
for $\tz=[z,\ze] \in \R^{n+1}$.

The {\it affine Weyl group} $W^a$ is generated by all $s_{\tal}$
(simple reflections $s_j=s_{\al_j} \for 0 \le j \le n$
are enough).
It is
the semi-direct product $W\lsmash A$, where the non-affine Weyl  
group $W$ is the span of $s_\al,
\al \in R_+$.
Here and futher we identify $b\in B$ with the corresponding
translations. For instance,
$$
\eqalignno{
& a =\ s_{\al}s_{[\al,1]}=\ s_{[-\al,1]}s_{\al}\for a=\al^{\vee},
\ \al\in R.
&\eqnu
\label\baction\eqnum*
}
$$

The {\it extended Weyl group} $ W^b$ generated by $W\and B$
 is isomorphic to $W\lsmash B$:
$$
\eqalignno{
&(wb)([z,\ze])\ =\ [w(z),\ze-(z,b)] \for w\in W, b\in B.
&\eqnu
}
$$

 Given $b_+\in B_+$, let
$$
\eqalignno{
&\om_{b_+} = w_0w^+_0  \in  W,\ \pi_{b_+} =
b_+(\om_{b_+})^{-1}
\ \in \ W^b, \ \om_i=\om_{b_i},\pi_i=\pi_{b_i},
&\eqnu
\label\wo\eqnum*
}
$$
where $w_0$ (respectively, $w^+_0$) is the longest element in $W$
(respectively, in $ W_{b_+}$ generated by $s_i$ preserving $b_+$) 
relative to the 
set of generators $\{s_i\}$ for $i >0$.

The elements $\pi_r\equal\pi_{b_r}, r \in O$ leave $\Ga^a$ invariant 
and form a group denoted by $\Pi$, 
 which is isomorphic to $B/A$ by the natural 
projection $\{b_r \to \pi_r\}$. As to $\{\om_r\}$,
they preserve the set $\{-\th,\al_i, i>0\}$.
The relations $\pi_r(\al_0)= \al_r= (\om_r)^{-1}(-\th) 
$ distinguish the
indices $r \in O^*$. Moreover (see e.g. [C2]):
$$
\eqalignno{
& W^b  = \Pi \lsmash W^a, \where
  \pi_rs_i\pi_r^{-1}  =  s_j \if \pi_r(\al_i)=\al_j,\  0\le j\le n.
&\eqnu
}
$$

Given $\nu\in\nu_R,\  r\in O^*,\  \tw \in W^a$, and a reduced  
decomposition $\tw\ =\ s_{j_l}...s_{j_2} s_{j_1} $ with respect to
$\{s_j, 0\le j\le n\}$, we call $l\ =\ l(\hw)$ the {\it length} of 
$\hw = \pi_r\tw \in W^b$. Setting
$$
\eqalign{
\la(\hw) = &\{ \tal^1=\al_{j_1},\
\tal^2=s_{j_1}(\al_{j_2}),\ 
\tal^3=s_{j_1}s_{j_2}(\al_{j_3}),\ldots \cr
&\ldots,\tal^l=\tw^{-1}s_{j_l}(\al_{j_l}) \},
}
\eqnu
$$
\label\tal\eqnum*
one can represent 
$$
\eqalign
{
&l=|\la(\hw)|=\sum_\nu l_\nu, \for l_\nu = l_\nu(\hw)=|\la_\nu(\hw)|,\cr
&\la_\nu(\hw) = \{\tal^{m},\ \nu(\tal^{m})= \nu(\tal_{j_m})= \nu\},
1\le m\le l,
}
\eqnu
\label\laset\eqnum*
$$
where $|\ |$  denotes the  number of elements,
 $\nu([\al,k]) \equal \nu_{\al}$.   

Let us introduce
the following {\it affine 
action} of $W^b$ on $z \in \R^n$:
$$
\eqalign{
& (wb)\langle z \rangle \ =\ w(b+z),\ w\in W, b\in B,\cr
& s_{\tal}\langle z\rangle\ =\ z - ((z,\al)+k)\al^\vee,
\ \tal=[\al,k]\in R^a,
}
 \eqnu
\label\afaction\eqnum*
$$
and the pairing $([z,\zeta], z'+d)
\equal (z,z')+\zeta$, where we treat $d$ formally (see e.g. [K]).
The connection with (\ref\saction,\ref\baction) is
as follows: 
$$
\eqalign{
& (\hw([z,\zeta]),\hw\langle z' \rangle+d) \ =\ 
([z,\zeta], z'+d) \for \hw\in W^b. 
}
 \eqnu
\label\dform\eqnum*
$$

Using the affine Weyl chamber 
$$
\eqalignno{
&C^a\ =\ \bigcap_{j=0}^n L_{\al_j},\ L_{\tal}=\{z\in \R^n,\ 
(z,\al)+k>0 \},
}
$$
$$
\eqalign{
\la_\nu(\hw)\ & =\ \{\tal\in R^a_+, \ \langle  C^a \rangle 
\not\subset \hw\langle L_{\tal}\rangle, \ \nu({\tal})=\nu \} \cr
& =\ \{\tal\in R^a_+, \ l_\nu( \hw s_{\tal}) < l_\nu(\hw) \}.
}
 \eqnu
\label\lambda\eqnum*
$$
It coincides with (\ref\laset) due to the relations 
$$
\eqalign{
&  \la_\nu(\hw\hu) = \la_\nu(\hu) \cup
\hu^{-1}(\la_\nu(\hw)),\ 
\la_\nu(\hw^{-1}) = -\hw(\la_\nu(\hw))\cr
&\if l_\nu(\hw\hu)=
l_\nu(\hw)+l_\nu(\hu).
} 
\eqnu
\label\ltutw\eqnum* 
$$
The following proposition is from [C4].

\proclaim {Proposition}
 Given $ b\in B$, the decomposition $b= \pi_b\om_b,
\om_b \in W$
can be uniquely determined  from the following equivalent conditions

{i)\ \  } 
$l(\pi_b)+l(\om_b)\ =\ l(b)$ and $l(\om_b)$ is the biggest possible,

{ii)\ }
$\om_b(b) = b_-\in B_-$ and $l(\om_b)$ is the smallest possible, 

{iii) }
$\pi_b\langle 0\rangle\ =\ b \and \la(\pi_b)\cap R\ =\ \emptyset$.
\endproclaim
\proofbox
\label\PIOM\theoremnum*

We will also use  that 
$$
\eqalignno{
\la(b) = \{ \tal,\  &( b, \al )>k\ge 0 \if \al\in R_+,
&\eqnu \cr
&( b, \al )\ge k> 0 \if \al\in R_- \},
\label\lambi\eqnum*\cr 
\la(\pi_b) = \{ \tal,\ \al\in R_-,\ &( b_-, \al )>k> 0 \if (\al,b)<0,
&\eqnu 
\label\lambpi\eqnum* \cr   
&( b_-, \al )\ge k > 0 \if (\al,b)>0 \}, \and \cr
\la(\pi_b^{-1}) = \{ \tal,\  -&(b,\al)>k\ge 0 \} \for \tal=[\al,k]\in R^a_+.
&\eqnu 
\label\lapimin\eqnum* 
}
$$

%\vfil
{\bf Convexity.} Let us introduce two orderings
on $B$. Here and further $b_{\pm}$ are the unique elements
from $B_{\pm}$ which belong to the orbit $W(b)$. Namely,
$b_-=\om_b\pi_b=\om_b(b)$, $b_+=w_0(b_-)= \om_{-b}(b).$
So the equality   $c_-=b_- $ (or $c_+=b_+ $) means that $b,c$
belong to the same orbit.
Set 
$$
\eqalignno{
&b \le c, c\ge b \for b, c\in B \iif c-b \in A_+,
&\eqnu 
\label\order\eqnum*
\cr
&b \preceq c, c\succeq b \iif b_-< c_- \hbox{\ \ or\ \ }
b_-=c_- \hbox{\ and\ } b\le c.
&\eqnu
\label\succ\eqnum*
}
$$
We  use $<,>,\prec, \succ$ respectively if $b \neq c$. 
For instance, 
$$c\succ b_+\Leftrightarrow b_+>W(c)>b_-,\ \ 
c\succeq b_-\Leftrightarrow c\in W(b_-) \hbox{\ or\ }  
c\succ b_+.
$$
The following sets  
$$
\eqalign{
&\si(b)\equal \{c\in B, c\succeq b\},\ 
\si_*(b)\equal \{c\in B, c\succ b\}, \cr  
&\si_+(b)\equal \{c\in B, c_->b_-\}\ =\ \si_*(b_+). 
}
\eqnu
\label\cones\eqnum*
$$
are convex. Moreover $\si_+$ is $W$-invariant. 
By {\it convex}, we mean that if
$ c, d= c+r\al^{\vee}\in \si$  
for $\al\in R_+, r\in \Z_+$, then
$$
\eqalignno{
&\{c,\ c+\al^{\vee},...,c+(r-1)\al^{\vee},\ d\}\subset \si.  
&\eqnu
\label\convex\eqnum*
}
$$
The elements from $\si(b)$
strictly between $c$ and $d$ (i.e.
$c+q\al,\ $ $0<q<r$) belong to $\si_+(b)$. 

\proclaim{Proposition }
a) Let $\hu= s_{\tal^{i_m}}...s_{\tal^{i_1}}\pi_b$, where $i_p$
are from
any sequence $1\le i_1<i_2<\ldots< i_m\le l=l(b)$ in
a reduced decomposition of $\hw=\pi_b^{-1}$
(see (\ref\tal)). In other words, 
$\hu$ is obtained by crossing out any number of $\{s_{j}\}$
from a reduced decomposition of $\pi_b$.
Then $c\equal\hu\langle 0\rangle \in \si_*(b)$.  Moreover, 
 $c\in \si_+(b)$
if and only if at least one of $\tal^{i_p}=[\al,k]$ 
for $1\le p\le m$ has $k>0$.  

b) If $c, b$ belong to the same $W$-orbit then the converse
is true. Namely, setting
$\om_{bc}\equal\pi_b\pi_c^{-1}$, 
the following relations are equivalent: 

(i)\ \  $c\succ b$ (which means that $c> b$),

(ii)\  $(\al,c)>0$ for all $\al\in \la(\om_{bc})$,

(iii) $l(\pi_b)\ =\ l(\om_{bc})+l(\pi_c)$,

It also results from (i) that
$\om_{bc}$ is the smallest possible element $w\in W$
such that $b=w(c)$. 
\label\BSTAL\theoremnum*
\endproclaim

{\it Proof.} Assertion a) is
a variant of Proposition 1.2 from [C4]. 
For the sake of completeness we
 will outline the proof of b).
Taking $u(c)\le b<c$, we will check (ii),(iii) by
induction supposing that
$\{u'(c)\le b'<c\}\Rightarrow$ \{(ii),(iii)\} 
 for all $b',u'$  such that
$l(u')<l(u)$,
which is obvious when $l(u')=0$.
  
Setting $\be=u(\al)$ for $\al\in \la(u)$, 
$u(s_{\al}(c))=u(c)-(\al,c)\be^\vee$ and
$\be\in R_-$ (see the definition of $\la(\al)$). One can
assume that $(\al,c)>0$ for all such $\al$. 
Otherwise $us_\al(c)\le u(c)\le c$
and we can argue by induction.
Applying (\ref\ltutw) and (\ref\lambi), 
we see that $l(uc)=l(u)+l(c)$.
Indeed, the intersection of  $\la(c)$ and
$$
c^{-1}(\la(u))=\{[\al,(c,\al)], \al\in \la(u)\} $$
is empty. Hence the product $u\pi_c$ is reduced (i.e.
$l(u\pi_c)=l(u)+l(\pi_c)$)  and 
$\la(u\pi_c)=\om_c c^{-1}( \la(u))\cup\la(\pi_c)$
contains no roots from $R_+$. Finally, Proposition \ref\PIOM
leads to (iii) (and the uniqueness of $u$ of the minimal
possible length). 
This reasoning gives the equivalence of (ii) and (iii) as well.
Assertion (i) readily results from (ii).
\proofbox

We will also use (cf. Proposition 5.2, [C4]) the relations 
$\pi_b=\pi_r\pi_c$ for $b=\pi_r\langle c\rangle$ and 
any $c\in B, r\in O$ and the
equivalence of the following three conditions: 
$$
\eqalignno{
& (\al_j,c+d)>0\Leftrightarrow \al_j\not\in \la(\pi_c^{-1})
\Leftrightarrow \{s_j\pi_c=\pi_b, c\succ b\}  
&\eqnu
\label\aljb\eqnum*
}
$$
for $0\le j\le n$.
When $j>0$ it is a particular case of Proposition \ref\BSTAL b).
Assuming that $(\al_0,c+d)=1-(\th,c)>0$,
$$
b=s_0\langle c\rangle = c+(\al_0,c+d)\th > c >c-\th=
s_\th (b).
$$
Hence $c\in \si_+(b)$. If the product $s_0\pi_c$
is reducible then we can apply statement a)
to come to a contradiction. Therefore
$s_0\pi_c=\pi_b$, since $s_0$ is simple.  
The remaining implications are obvious.

%
%
%		
%                  Section 2
%
%
%\vskip 10pt
\section{ Intertwining operators} 
We put    
$m=2 \for D_{2k} \and C_{2k+1},\ m=1 \for C_{2k}, B_{k}$,
otherwise $m=|\Pi|$. Let us set
$$
\eqalignno{
&   t_{\tal} = t_{\nu(\tal)},\ t_j = t_{\al_j},
\where \tal \in R^a,\ 0\le j\le n, \cr 
& X_{\tb}\ =\ \prod_{i=1}^nX_i^{k_i} q^{ k} 
\if \tb=[b,k],
&\eqnu \cr 
&\for b=\sum_{i=1}^nk_i b_i\in B,\ k \in {1\over 2m}\Z.
}
$$
\label\Xde\eqnum*
Here and futher
$ q,\{ t_\nu , \nu \in \nu_R \},$ $X_1,\ldots,X_n$
are considered as  independent variables,
$ \C_{ q,t}$
 is the field of rational 
functions in terms of $ q^{1/2m},
\{t_\nu^{1/2} \}$,
$\C_{q,t}[X] = \C_{q,t}[X_b]$  means the algebra of
polynomials in terms of $X_i^{\pm 1}$ 
with the coefficients from $ \C_{ q,t}$.

We will keep the notations:
 $$
([a,k],[b,l])=(a,b),\
[\al,k]^\vee=2[\al,k]/(\al,\al), \ \nu_{\al^\vee}=\nu_\al,
$$
$a_0=\al_0, X_0=X_{a_0}$, and use the involution 
$$
O^*\ni r\to r^*,\ \al_{r^*} \equal \pi_r^{-1}(\al_0).
$$ 
Check that $\om_r\om_{r^*}=1=\pi_r\pi_{r^*}$.

\proclaim{Definition }
 The  double  affine Hecke algebra $\HH\ $
(see [C1,C2]) 
is generated over the field $ \C_{ q,t}$ by 
the elements $\{ T_j,\ 0\le j\le n\}$, 
pairwise commutative $\{X_b, \ b\in B\}$ satisfying (\ref\Xde),
 and the group $\Pi$ where the following relations are imposed:

(o)\ \  $ (T_j-t_j^{1/2})(T_j+t_j^{-1/2})\ =\ 0,\ 0\ \le\ j\ \le\ n$;

(i)\ \ \ $ T_iT_jT_i...\ =\ T_jT_iT_j...,\ m_{ij}$ factors on each side;

(ii)\ \   $ \pi_rT_i\pi_r^{-1}\ =\ T_j \if \pi_r(\al_i)=\al_j$; 

(iii)\  $T_iX_b T_i\ =\ X_b X_{a_i}^{-1} \if (b,\al_i)=1,\
1 \le i\le  n$;

(iv)\  $T_0X_b T_0\ =\ X_{s_0(b)}\ =\ X_b X_{\th} q^{-1}
\if (b,\th)=-1$;

(v)\ \ $T_iX_b\ =\ X_b T_i$ if $(b,\al_i)=0 \for 0 \le i\le  n$;

(vi)\ $\pi_rX_b \pi_r^{-1}\ =\ X_{\pi_r(b)}\ =\ X_{\om^{-1}_r(b)}
 q^{(b_{r^*},b)},\  r\in O^*$.
\label\double\theoremnum*
\endproclaim
\proofbox

Given $\tw \in W^a, r\in O,\ $ the product
$$
\eqalignno{
&T_{\pi_r\tw}\equal \pi_r\prod_{k=1}^l T_{i_k},\where 
\tw=\prod_{k=1}^l s_{i_k},
l=l(\tw),
&\eqnu
\label\Tw\eqnum*
}
$$
does not depend on the choice of the reduced decomposition
(because $\{T\}$ satisfy the same ``braid'' relations as $\{s\}$ do).
Moreover,
$$
\eqalignno{
&T_{\hv}T_{\hw}\ =\ T_{\hv\hw}\  \hbox{ whenever}\ 
 l(\hv\hw)=l(\hv)+l(\hw) \for
\hv,\hw \in W^b.
 %&(2.7)}
&\eqnu}
$$
\label\TT\eqnum*
  In particular, we arrive at the pairwise 
commutative elements 
$$
\eqalignno{
& Y_{b}\ =\  \prod_{i=1}^nY_i^{k_i} \if  
b=\sum_{i=1}^nk_ib_i\in B,\where  
 Y_i\equal T_{b_i},
&\eqnu
\label\Yb\eqnum*
}
$$
satisfying the relations
$$
\eqalign{
&T^{-1}_iY_b T^{-1}_i\ =\ Y_b Y_{a_i}^{-1} \if (b,\al_i)=1,
\cr
& T_iY_b\ =\ Y_b T_i \if (b,\al_i)=0, \ 1 \le i\le  n.}
%\eqno(2.9)
\eqnu
$$

\comment
Let us introduce the following elements from
$\C_t^n$:
$$
\eqalign{
&t^{\pm\rho}\equal (l_t(b_1)^{\pm 1},\ldots,l_t(b_n)^{\pm 1}),\where\cr
&l_t(\hw)\equal \ \prod_{\nu\in\nu_R} t_\nu^{l_\nu(\hw)/2},\
\hw\in W^b,
}
\eqnu
\label\qlen\eqnum*
$$
and the corresponding {\it evaluation maps}:
$$
\eqalign{
&X_i(t^{\pm\rho})= l_t(b_i)^{\pm 1} = Y_i(t^{\pm\rho}),\ 1\le i\le n.
}
\eqnu
\label\eval\eqnum*
$$
For instance, $X_{a_i}(t^{\rho})\ =\ l_t(a_i)= t_i$ (see (\ref\lb)).

We will establish the  duality of non-symmetric polynomials applying
the following theorem ([C2],[C3]).
\proclaim {Theorem}
i) The elements $H \in \HH\ $  have
the unique decompositions
$$
\eqalignno{
&H =\sum_{w\in W }  g_{w}  T_{w} f_w,\ 
g_{w} \in \C_{ q,t}[X],\ f_{w} \in \C_{ q,t}[Y].  
&\eqnu
}
$$

ii) The   map 
$$
\eqalign{
 \vph: &X_i \to Y_i^{-1},\ \  Y_i \to X_i^{-1},\  \ T_i \to T_i, \cr
&t_\nu \to t_\nu,\
 q\to  q,\ \nu\in \nu_R,\ 1\le i\le n.
}
\eqnu
\label\vph\eqnum*
$$
can be extended to an anti-involution 
($\vph(AB)=\vph(B)\vph(A)$)
 of \HH\ .

iii) The  linear functional on \HH\ 
$$
\eqalignno{
&[\![ \sum_{w\in W }  g_{w}  T_{w} f_{w}]\!]\ =\
\sum_{w\in W} g_{w}(t^{-\rho}) l_t(w) f_{w}(t^{\rho})
&\eqnu
\label\brack\eqnum*
}
$$
is invariant with respect to $\vph$. The bilinear form 
$$
\eqalignno{
&[\![ G,H]\!]\equal [\![ \vph(G)H]\!],\
G,H\in \HH\ ,
&\eqnu
\label\form\eqnum*
}
$$
is symmetric ($[\![ G,H]\!]= [\![ H,G]\!]$)
and non-degenerate.
\endproclaim
\label\dual\theoremnum*
\proofbox

The map $\vph$ is the composition of the 
involution (see [C1])
\endcomment

The following maps can be extended to involutions of
\HH\ (see [C1,C3]):
$$
\eqalignno{
  \vep:\ &X_i \to Y_i,\ \  Y_i \to X_i,\  \ T_i \to T_i^{-1},
&\eqnu
\label\vep\eqnum*
 \cr
&t_\nu \to t_\nu^{-1},\
 q\to  q^{-1},\cr
 \tau: \ &X_b \to X_b,\ \ Y_r \to X_rY_r q^{-(b_r,b_r)/2},\ 
Y_\th \to X_0^{-1}T_0^{-2}Y_\th, &\eqnu
\cr
&T_i\to T_i,\ t_\nu \to t_\nu,\
 q\to  q, \cr
&1\le i\le n,\ r\in O^*,\ X_0=qX_\th^{-1}.
\label\tau\eqnum*
}
$$

Let us give some explicit formulas:
$$
\eqalign{
&\vep(T_0)\ =\ X_\th T_0^{-1}Y_\th\ =\ X_\th T_{s_\th},\ 
\vep(\pi_r)\ =\ X_rT_{\om_r^{-1}},\cr
&\tau(T_0)\ =\ X_0^{-1}T_0^{-1},\ \tau(\pi_r)\ =
\ q^{-(b_r,b_r)/2}X_r\pi_r= q^{(b_r,b_r)/2}\pi_rX_{r^*}^{-1},\cr
&\pi_r X_{r^*}\pi_r^{-1}\ =\ q^{(b_r,b_r)}X_{r}^{-1},\ 
X_{r^*}T_{\om_r}X_r\ =\ T_{\om_{r^*}}^{-1}.
}
\eqnu
\label\vphto\eqnum*
$$

Theorem 2.3 from [C3] says that the map
$$
\eqalign{
& \Bigl(\matrix  0 &-1\\ -1& 0\endmatrix \Bigr) \to \vep,\ 
  \Bigl(\matrix 1& 1\\  0& 1 \endmatrix \Bigr) \to \tau
}
\eqnu
\label\glz\eqnum*
$$
can be extended to a homomorphism of $GL_2(\Z)$ up to conjugations 
by the central elements from the group generated
by $T_1,\ldots,T_n$. 

The  involution $\eta=\tau^{-1}\vep\tau$
corresponding to the
matrix $\Bigl(\matrix -1&  0\\  1& 1 \endmatrix \Bigr)$
will play an important role in the paper:
$$
\eqalignno{
  \eta: \ &X_r\to  q^{(b_r,b_r)/2}X_r^{-1}Y_r\ =\ \pi_rX_{r^*}T_{\om_r},
&\eqnu
\cr
 &Y_r \to q^{(b_r,b_r)/2}X_r^{-1}Y_rX_r\  
 =\ \pi_rT_{\om_{r^*}}^{-1},\cr 
&Y_\th \to T_0^{-1}T_{s_\th}^{-1},\ 
T_j\to T_j^{-1} (0\le j\le n),\cr
&\pi_r\to\pi_r (r\in O^*),\ t_\nu \to t_\nu^{-1},\
 q\to  q^{-1}.
\label\eta\eqnum*
}
$$

We note that $\vep$ and $\eta$ commute with 
the main anti-involution $^*$ from [C2]:
$$
\eqalign{
  & X_i^*\ =\  X_i^{-1},\ \  Y_i^*\ =\  Y_i^{-1},\  \
 T_i^* \ =\  T_i^{-1}, \cr
&t_\nu \to t_\nu^{-1},\
 q\to  q^{-1},\ 0\le i\le n,\ 
(AB)^*=B^*A^*.
}
\eqnu
\label\star\eqnum*
$$

The {\it $X$-intertwiners} (see e.g. [C2,C5,C6])
are introduced as follows:
$$
\eqalignno{
&\Phi_j\ =\ 
T_j + (t_j^{1/2}-t_j^{-1/2})(X_{a_j}-1)^{-1},
\cr
& G_j=\Phi_j (\phi_j)^{-1},\
\tG_j=(\phi_j)^{-1}\Phi_j,\ 
 &\eqnu
\cr 
&\phi_j\ =\  t_j^{1/2} + 
(t_j^{1/2} -t_j^{-1/2})(X_{a_j}-1)^{-1},
\label\Phi\eqnum*
}
$$
for $0\le j\le n$. They belong to the 
extension of \HH\ by the field $\C_{q,t}(X)$
of rational functions in $\{X\}$. The elements
 $G_j$ and $G'_j$
satisfy the same relations
as $\{s_j,\pi_r\}$ do, $\{\Phi_j\}$ satisfy the relations 
for $\{T_j\}$ (i.e. the  homogeneous Coxeter relations and 
those with $\pi_r$).  Hence the elements 
$$
G_{\hw}\ =\ \pi_r G_{j_l}\cdots G_{j_1},
\where \hw=\pi_r s_{j_l}\cdots s_{j_1}\in W^b,
\eqnu
\label\Phiprod\eqnum*
$$
are well-defined and $G$ is a homomorphism of $W^b$.
The same holds for $\tG$. As to $\Phi$, the decomposition of $\hw$
should be reduced.

The simplest way to see this is to use the following 
property of $\{\Phi\}$ which fixes them uniquely up to left
or right multiplications by functions of $X$:
$$
\eqalignno{
&\Phi_{\hw} X_b\ =\ X_{\hw(b)}\Phi_{\hw},\ \hw\in W^b.
&\eqnu 
\label\Phix\eqnum*
}
$$
One first checks (\ref\Phix) for $s_j$ and $\pi_r$, then
observes that $\Phi$ from (\ref\Phiprod) satisfy (\ref\Phix) for
any choice of the reduced decomposition, and uses the 
normalizing conditions to see that they are uniquely
determined from the intertwining relations  (\ref\Phix).

We note that $\Phi_j,\phi_j$ are self-adjoint with respect to 
the anti-involution (\ref\star). Hence
$$
\eqalignno{
& \Phi_{\hw}^* = \Phi_{\hw^{-1}},\ 
G_{\hw}^*\ =\ \tG_{\hw^{-1}},  \ \hw\in W^b.
&\eqnu 
\label\Phistar\eqnum*
}
$$
It follows from the quadratic relations for $T$.

To define the {\it $Y$-intertwiners} we apply
the involution $\vep$ to $\Phi_{\hw}$ and to $G,\tG$.
The formulas can be easily calculated using
(\ref\vphto). In the case of $GL_n$ one gets 
the intertwiners from [Kn].
For $w\in W$, we just need to replace $X_b$
by $Y_b^{-1}$ and conjugate $q,t$ (cf. [C4]). 
However it will be  more convenient to consider 
$\eta(\Phi)$ instead
of $\vep(\Phi)$ to create the Macdonald polynomials.
Both constructions gives the intertwiners satisfying
the $\ast$-relations from (\ref\Phistar).

%
%
%		Section 3
%
%
%\vskip 10pt
\section { Standard representations} 
It was observed in [C4], Section 5 that there is a natural passage from 
the representation of \HH\ in polynomials to a representation
in functions on $W^b$. We will continue this line, beginning with
the construction of the {\it basic representaions} of level $0,1$.
Setting 
$$
\eqalignno{
& x_{\tb}=  \prod_{i=1}^nx_i^{k_i} q^{ k} \if 
\tb=[b,k],\
b=\sum_{i=1}^nk_i b_i\in B,\ k \in {1\over 2m}\Z,
&\eqnu 
\label\xde\eqnum*}
$$
for independent $x_1,\ldots,x_n$, we
 consider $\{X\}$ as  operators acting in $\C_{q,t}[x]=$
$\C_{q,t}[x_1^{\pm 1},$ $\ldots,x_n^{\pm 1}]$:
$$
\eqalignno{
& X_{\tb} (p(x))\ =\ x_{\tb} p(x),\    p(x) \in
\C_{q,t} [x]. 
&\eqnu}
$$
\label\X\eqnum*
The elements $\hw \in W^b$ act in $\C_{ q}[x]$
in two ways:
$$
\eqalignno{
&\hw(x_{\tb})\ =\ x_{\hw(\tb)},\ \hw\llan x_{\tb}\rran\ =\ 
x_{\hw\llan \tb \rran}, \where\cr
& s_{\tal}\llan \tb\rran= \tb-(\tal,b+d)\tal^\vee,
\ 
a \llan \tb\rran= a(\tb)+[a,-(a,a)/2]
&\eqnu
\label\afxact\eqnum*
}
$$
for $a\in B,\tal\in R^a$.
More generally, we can replace in (\ref\afxact) 
$d$ by $ld$ and
$[,]$ by $l[,]$ 
(the action of level $l$ like for Kac-Moody algebras) 
but only $l=0,1$
will be used in this paper. The most general action
depends on an element of $SL(2,\Z)$.

Thus (\ref\afxact) is an extension of the affine action
$\langle\ \rangle$
from (\ref\afaction)  to $R^{n+1}\ni [b,k]$.
The affine action on functions will always
mean (\ref\afxact).
In particular,
$$
\eqalignno{
&\pi_r(x_{b})\ =\ x_{\pi_r({b})}\ =\ 
x_{\om^{-1}_r(b)} q^{(b_{r^*},b)}\cr 
&\pi_r\llan x_{b}\rran\ =\ x_{\pi_r({b})}x_{b_r}q^{-(b_r,b_r)/2}, 
&\eqnu
\label\pi\eqnum*
}
$$
where (we remind) $\al_{r^*}= \pi_r^{-1}(\al_0), \ r\in O^*$.

Respectively, the {\it Demazure-Lusztig operators} (see  
[C2]) 
$$
\eqalignno{
&\hT_j\  = \  t_j ^{1/2} s_j\ +\ 
(t_j^{1/2}-t_j^{-1/2})(X_{a_j}-1)^{-1}(s_j-1),
\ 0\le j\le n.
&\eqnu
\label\Demaz\eqnum*
}
$$
act in $\C_{ q,t}[x]$ according to the level $l=0,1$.

We note that only $\hT_0$ depends on $ q$: 
$$
\eqalign{
&\hT_0\  =\  t_0^{1/2}s_0\ +\ (t_0^{1/2}-t_0^{-1/2})
( q X_{\th}^{-1} -1)^{-1}(s_0-1),\cr
&\where
s_0(x_b) = x_b x_{\th}^{-(b,\th)} q^{(b,\th)},\
s_0\llan x_b\rran = 
x_b (qx_{\th}^{-1})^{(b,\th)-1}. 
}
\eqnu
$$
 
\proclaim{Theorem }
 The map $ T_j\to \hT_j,\ X_b \to X_b$ (see (\ref\Xde,\ref\X)),
$\pi_r\to \pi_r$  (see (\ref\pi)) gives two
\HH\ - modules 
$\v_0\simeq \C_{ q,t}[x]\simeq \v_1$ (for\ $l=0,1$) 
over $ \C_{ q,t}$. 
The action of $H\in \HH\ $ in $\v_1$ coincides
with the action of $\tau(H)$ in $\v_0$.
Generally speaking, one can introduce the module $\v_g$
for $g\in SL(2,\Z)$ acting on $H$ by the 
outer automorphism corresponding to $g$ (see (\ref\glz)
above and Theorem 4.3, [C3]).  
These representations are faithful and 
remain faithful when   $  q,t$ take  any nonzero
values assuming that
 $ q$ is not a root of unity (see [C2]).
The representation $\v_0$ is induced from the character
$\{T_j\to t_j, \pi_r\to 1\}$. Namely, the image $\hat{H}$
is uniquely determined from the following condition:
$$
\eqalign{
&\hat{H}(f(x))\ =\ g(x)\for H\in \HH\ \if Hf(X)-g(X)\ \in\cr
 & \sum_{i=0}^n \HH\ (T_i-t_i)+
\sum_{r\in O^*} \HH\ (\pi_r-1).}
\eqnu
\label\hat\eqnum*
$$ 
\endproclaim
\proofbox
\label\faith\theoremnum*

To make the statement about $\v_1$ quite obvious let us 
introduce the {\it Gaussian} 
$\ga\ =\ \hbox{Const}\ q^{\Sigma_{i=1}^n z_i z_{\al_i}/2}$,
where formally
$$
x_b= q^{z_b}, \ z_{a+b}=z_a+z_b, \ z_i=z_{b_i},\ 
(wa)(z_b)= z_{w(b)}-(a,b), \ a,b\in \R^n.
$$
More exactly, it is a $W$-invariant
solution of  the following  difference equations:
$$
\eqalign{
&b_j(\ga)\ =\ \hbox{Const}\  q^{(1/2)\Sigma_{i=1}^n (z_i-(b_j,b_i))
(z_{\al_i}- \de_i^j)}\ =\cr
&q^{-z_j+ (b_j,b_j)/2 }\ga \ =\  q^{(b_j,b_j)/2}
x_j^{-1}\ga  \for 1\le j\le n.   
}
\eqnu
\label\gauss\eqnum*
$$

The Gaussian  commutes with $T_j \for 1\le j\le n$
because it is $W$-invariant.
A straightforward calculation  gives that 
$$
\eqalign{
&\ga(X)T_0\ga(X)^{-1}\ =\ X_0^{-1}T_0^{-1}\ =\  \tau(T_0) ,\cr
&\ga(X)Y_r\ga(X)^{-1}\ =\  q^{-(b_r,b_r)/2}  X_rY_r \ =
\ \tau(Y_r), r\in O.
}
\eqnu
\label\gato\eqnum*
$$
Hence the conjugation by $\ga$ induces $\tau$.
We can put in the following way. There is a formal 
$\HH$-homomorphism:
$$
\eqalign{
&\v_0\ni v\to\hv\equal v\ga^{-1}\in \v_1.
}
\eqnu
\label\vovone\eqnum*
$$
One has to complete $\v_{0,1}$ to make this map 
well-defined (see the discrete representations below).

We will later need an extended version of Proposition 3.6
from [C2].

\proclaim{Proposition }
a)The operators 
$\{ Y_i, 1\le i\le n\}$ acting in $\v_0$
 preserve $\Si(b)\equal
\oplus_{c\in \si(b)}\C_{q,t} x_{c}$ and
  $\Si_*(b)$ (defined for $\si_*(b)$)
for arbitrary $b\in B$. 

b)The operators
$\{T_j, 0\le j\le n\}$ acting in $\v_0$
preserve $\Si_+(b)=\Si_*(b_+)$:
$$
\eqalign{
\hT_j(x_b) &\hbox{\ mod\ } \Si_+(b)\ =\
t_j^{1/2}x_b\if (b,\al_j)=0, \cr 
&=\ t_j^{1/2}s_j(x_b)+(t_j^{1/2}-t_j^{-1/2})x_b
 \if (b,\al_j)<0,\cr
&=\ t_j^{-1/2}s_j(x_b) \if (b,\al_j)>0.
}
\eqnu
\label\tonx\eqnum*
$$

c) Coming to $\v_1$, if $(\al_j,b+d)>0\ (0\le j\le n)$
then 
 $$
\eqalign{
\hT_j\llan x_b\rran &\hbox{\ mod\ } \Si_+(s_j\lan b\ran)\ = \ 
t_j^{-1/2} s_j\llan x_b \rran.
}
\eqnu
\label\tmonx\eqnum*
$$
Otherwise, $(\al_j,b+d)\le 0$ and
$$
\eqalign{
&\hT_j\llan x_b\rran \in \Si( b)\for (\al_j,b+d)\le 0,\cr
&\hT_j\llan x_b\rran \ =\ t_j^{1/2}x_b\if (\al_j,b+d)= 0.
}
\eqnu
\label\tlonx\eqnum*
$$
\endproclaim
\label\TONX\theoremnum*
{\it Proof.} Due to Proposition 3.3 from [C4]
it suffices  to check  c) for $j=0$.
The first inequality, the definition of
$\hT_0\llan x_b\rran=\sum_{c\in B} u_{bc}x_c$, 
and (\ref\aljb) readily
give that (for nonzero $u$) $c=b+r\th\ (r\in \Z)$ and
 $$
\eqalign{
s_\th(b')=b-\th< b\  < c\ \le\ b+(\al_0,b+d)\th=
s_0\lan b\ran\equal b'.
}
\eqnu
\label\soineq\eqnum*
$$
Hence $c\in \si_+(b')$
if $c\neq b'$. The coefficient $u_{bb'}$ equals $t_0^{-1/2}$. 
Let $(\al_0,b+d)\le 0$. Then 
$$
\eqalign{
s_\th(b)=b-(b,\th)\th \  < c\ \le\ b \and c\in \{\Si_+(b)\cup b\}
\in\Si(b)
}
\eqnu
\label\soineqq\eqnum*
$$
(cf. Proposition \ref\BSTAL, a)).
\proofbox

%\vfil
{\bf Discretization.} We go to the lattice
version of the  functions and operators. Let $\xi$ be 
a "generic" character of $\C[x]$:
$$
x_a(\xi)\equal \prod_{i=1}^n \xi_i^{k_i}, \ a=\sum_{i=1}^n k_i b_i\in B,
$$
for independent parameters
$\xi_i$. The 
 discretizations of 
functions $g(x)$ in $x\in \C^n$  and the operators
from the algebra $\a\equal \oplus_{\hu\in W^b}\C_{q,t}(X)\hu $,
are described by the formulas: 
$$
\eqalignno{
&{}^\de x_a(bw)=x_a(q^b w(\xi))=
q^{(a,b)}x_{w^{-1}(a)}(\xi),\cr 
 &({}^\de \hu ({}^\de g))(bw)\ =\ {}^\de g(\hu^{-1}bw).
&\eqnu
\label\deltag\eqnum*
}
$$
For instance,
$({}^\de X_a({}^\de g))(bw)\ =\  x_{a}( bw)
\  g(bw)$ (we will sometimes omit
${}^\de$ and put $g(\hw)$ instead of  ${}^\de g(\hw)$).

The image of $g\in\C_{q,t}(x)$  belongs to 
the space $\f_\xi\equal \hbox{Funct}(W^b,\C_\xi)$
of $\C_\xi$-valued functions on $W^b$, where
$\C_\xi\equal\C_{q,t}(\xi_1,\ldots,\xi_n)$.

Considering the discretizations of 
 operators
$\hH$ for $H\in \HH$ we come to the {\it functional representation} 
of $\HH$ in $\f_\xi$.  

Similarly, introducing the group algebra 
$\C_\xi[W^b]=\oplus_{\hw\in W^b}\C_\xi\de_{\hw}$ 
for (formal) {\it delta-functions},
we can consider the dual anti-action on the indices:
$$
\eqalign{
& {}_\de (g(x)\hu)(\sum_{{\hw}\in W^b}c_{\hw}\de_{\hw})\ =
\ \sum_{{\hw}\in W^b}c_{\hw} g( \hw)\de_{\hu^{-1}\hw},\ 
c_{\hw}\in \C_\xi.
}
\eqnu
\label\deltsuf\eqnum*
$$
Composing it with the  anti-involution of \HH\
$$
\eqalign{
& T_j^\diamond =T_j (0\le j\le n),\ \pi_r^\diamond
 = \pi_r^{-1} (r\in O),\
X_i^\diamond=X_i (0\le i\le n),
}
\eqnu
\label\antivee\eqnum*
$$
sending $q,t$ to $q,t$ (and $AB$ to $B^\diamond A^\diamond$), we
get the {\it delta-representation}  $\De_\xi$
of $\HH\ $ in $\C_\xi[W^b]$:
$$
\eqalign{
H\ \to\  _\de (\hat{H}^\diamond)\equal  \de(H) \for H\in \HH. 
}
\eqnu
\label\derep\eqnum*
$$

\comment
We will mostly use the discretizations 
$ \ep_b({\hw})= e_b ({\hw})/e_b (0)$
of the
renormalized Macdonald polynomials
 $\ep_b(x)= e_b(x)/e_b(t^{-\rho})$, and
especially
$\ep_b(\#c)=\ep_c(\#b),$ where
$\#c\equal\pi_c=c\om_c^{-1}.$
See (\ref\Yone),(\ref\rhob),
and Theorem \ref\DUAL.
Sometimes we drop $\#$ and write $\ep_b(c)$ instead of
$\ep_b(\#c)$.
Vice versa, we will consider  the sufficies $b$ as
 elements from $W^b$
via the same  map $b\to \#b$.
\endcomment

Explicitly, $^\de\pi_r=\pi_r=\de (\pi_r) \ , r\in O$, and for $\hw=bw$
$$
\eqalignno{
^{\de}(T_i(g))(\hw))\ =\ 
{t_i^{1/2}x_{a_i}(w(\xi))q^{(a_i,b)} - t_i^{-1/2}\over
{x_{a_i}(w(\xi))q^{(a_i,b)}- 1}  }\
&g(s_i \hw)\cr
-{t_i^{1/2}-t_i^{-1/2}\over
{x_{a_i}(w(\xi))q^{(a_i,b)}- 1}  }\
&g(\hw) \for 0\le i\le n,
&\eqnu \cr
\label\Tfunct\eqnum*
\de(T_i)(\de_{\hw})\ =\ 
{t_i^{1/2}x_{a_i}(w(\xi))q^{(a_i,b)} - t^{-1/2}\over
{x_{a_i}(w(\xi))q^{(a_i,b)}- 1}  }\
&\de_{s_i\hw}\cr
-{t_i^{1/2}-t_i^{-1/2}\over
{x_{a_i}(w(\xi))q^{(a_i,b)}- 1}  }\
&\de_{\hw} \for 0\le i\le n.
&\eqnu
\label\Tdelta\eqnum*
}
$$

There is a natural $\C_\xi$-linear
pairing between $\f_\xi$ and $\De_\xi$.
Given $g\in \hbox{Funct}(W^b,\C_\xi), \hw\in W^b$,  
$$
\eqalignno{
&\{g,\de_{\hw}\}\equal g(\hw),\ \{H(g),\de_{\hw}\}\ =\
\{g, H^\diamond(\de_{\hw})\},\ H\in \HH\ .
&\eqnu 
\label\depair\eqnum*
}
$$
It also gives a nondegenerate
pairing between $\v_{0}$ and $\De_\xi$.
For arbitrary operators $A\in \a$, the relation is
as follows: $ \{{}^{\de}A(g),\,\de_{\hw}\}\ =\
\{g,\, {}_{\de}A(\de_{\hw})\}$.

Let us extend the discretization map 
and the pairing with $\De_\xi$
 to  $\v_1$. We use the map from (\ref\vovone) for
 the $\de$-{\it Gaussian}:
$$
\eqalign{
&{}^\de \ga (bw)\equal q^{(b,b)/2}x_b(w(\xi)),
 }
\eqnu
\label\degauss\eqnum*
$$
which satisfies (\ref\gauss) and is a discretization of 
$\ga$ for a proper constant (cf. [C4], (6.20)).

The representations $\f_\xi$ and $\De_\xi$
can be introduced when $q,t,\{\xi_i\}$ are considered as
complex numbers ensuring that $x_{\ta}(\xi)\neq 1$ for
all $\ta\in (R^a)^\vee$.
Following Proposition 5.2 from [C4], let us specialize the definition
of $\De$ for $\xi=t^{-\rho}$. In this case
$$
\eqalignno{
&x_a(bw)=x_a(q^b t^{-w(\rho)})=
q^{(a,b)}\prod_\nu t_\nu^{-(w(\rho_\nu),a)}.
&\eqnu
\label\deltaf\eqnum*
}
$$
\proclaim{ Proposition} The $\HH$-module $\De(-\rho)\equal
\De_{t^{-\rho}}$ contains the
$\HH$ -sub\-mo\-dule $\De_\# \equal$ $ \oplus_{b\in B} \C\de_{\pi_b}$.
This also holds for any $q\in\C^*$ and generic $t$. Moreover,
$\De_\#$ is irreducible if and only if $q$ is not a root of unity.
\label\DELTAO\theoremnum*
\endproclaim
\proofbox 

%\vfill
When $q\to 0$ and $t$ is a power of prime $p$
the  action of the algebra 
$\h^a$ generated by $\{T_j, 0\le j \le n\}$ in $\De(-\rho)$
 coincides with
the standard action of the $p$-adic Hecke algebra 
$H(G/\!/I) \cong \h^a$
on  the (linear span of)  delta-functions on $I\backslash G/I \cong  W^a$.
Here $I$ is the Iwahori subgroup of the split semisimple $p$-adic
group $G$  (see [IM]).
However $\De_\#$ does not remain a submodule in this limit.

Multiplying the delta-functions on the right by the 
operator of $t$-sym\-met\-ri\-za\-tion we can get an $\h^a$-submodule
isomorphic to $\De_\#$ (upon the restriction to $\h^a$).
Its limit readily exists and coincides with 
the space of
delta-functions on $I\backslash G/K$ for
the maximal parahoric
subgroup $K$. However the latter space can be identified
with neither  spaces of delta-functions for smaller  subsets of $W^a$ 
(as in Proposition 
\ref\DELTAO). It is possible only for the $q$-deformation under
consideration. Practically, when  calculating with 
right $K$-invariant functions in the $p$-adic case one
 needs to consider their values on the whole $W^a$ 
(that is an obvious flaw since much fewer number of points 
is enough to reconstruct them uniquely).

\comment
Here we do not distinguish $\De_\#$ and
the representation of $\h^a$ in finitely supported
functions from $\f$. They are isomorphic (see below).

Let us introduce (formally) the element $d$ such that
$([z,k],d)=k$ for all $z\in \R^n$. Then
explicitly
$$
\eqalign{
&_{\de}T_i(\phi_{\#b})\ =\ 
{t_i^{1/2}t^{-(a_i,\om_b^{-1}(\rho))}q^{(a_i,b+d)} - t^{-1/2}\over
{t^{-(a_i,\om_b^{-1}(\rho))}q^{(a_i,b+d)}- 1}  }
\phi_{s_i\#b}\ -\cr
&{t_i^{1/2}-t_i^{-1/2}\over
{t^{-(a_i,\om_b^{-1}(\rho))}q^{(a_i,b+d)}- 1}  }
\phi_{\#b}, \for 0\le i\le n.
}
\eqnu
\label\Tdelta\eqnum*
$$
\endcomment

%
%
%		Section 4
%
%
%\vskip 10pt
\section { Orthogonality} 
The coefficient of $x^0=1$ ({\it the constant term})
of a polynomilal
$f\in  \C_{q,t}[x]$
will be denoted by $\langle  f \rangle_0$. Let
$$
\eqalign{
&\mu\ =\ \prod_{a \in R_+^\vee}
\prod_{i=0}^\infty {(1-x_aq_a^{i}) (1-x_a^{-1}q_a^{i+1})
\over
(1-x_a t_aq_a^{i}) (1-x_a^{-1}t_a^{}q_a^{i+1})},
}
\eqnu
\label\mu\eqnum*
$$
where $q_a=q_{\nu}=q^{2/\nu} \for \nu=\nu_a$.
 
The coefficients of $\mu_0\equal \mu/\langle \mu \rangle_0$
are from $\C(q,t)$, where the formula for the
constant term of $\mu$ is as follows
(see [C2]):
$$
\eqalign{
&\langle\mu\rangle_0\ =\ \prod_{a \in R_+^\vee}
\prod_{i=1}^\infty {(1-x_a(t^\rho)q_a^{i})^2
\over
(1-x_a(t^\rho) t_aq_a^{i}) (1-x_a(t^\rho) t_a^{-1}q_a^{i})}.
}
\eqnu
\label\consterm\eqnum*
$$
Here  
 $x_{b}(t^{\pm\rho}q^{c})=  
q^{(b,c)}\prod_\nu t_\nu^{\pm(b,\rho_\nu)}$. 

We note that
$\mu_0^*\ =\ \mu_0$  with respect to the involution 
$$
 x_b^*\ =\  x_{-b},\ t^*\ =\ t^{-1},\ q^*\ =\ q^{-1}.
$$
 
Setting 
$$
\eqalignno{
&\langle f,g\rangle_0\ =\langle \mu_0 f\ {g}^*\rangle_0\ =\ 
\langle g,f\rangle_0^* \for
f,g \in \C(q,t)[x], 
&\eqnu
\label\innerpro\eqnum*  
}
$$
 we  introduce the {\it non-symmetric Macdonald 
polynomials} $e_b(x),\   b \in B_-$, by means of 
the conditions
$$
\eqalignno{
&e_b-x_b\ \in\ \Si_*(b),\
\langle e_b, x_{c}\rangle = 0 \for c\in \si_* =
\{c\in B, c\succ b\} 
&\eqnu
\label\macd\eqnum*
}
$$
in the setup of Section 1.
They can be determined by the Gram - Schmidt process 
because the  pairing   
is non-degenerate 
 and form a 
basis in $\C(q,t)[x]$. 

This definition is due to Macdonald [M3] (for
 $t_\nu=q^k,\ k\in \Z_+$)
 who extended 
 Opdam's nonsymmetric polynomials introduced
in the degenerate (differential) case in [O2]. He also 
established the connection with the $Y$-operators.
The general case was considered in [C4].

 The notations are
from  Proposition \ref\PIOM and (\ref\rhor).
 We use the involution 
$ \bar{x}_a=x_a^{-1},\
\bar{q}= q,\ \bar{t}=t,\ a\in B$.

\proclaim {Proposition}
a) For any $H\in \HH\ $ and the anti-involution $^*$
from (\ref\star), $\langle \hH(f),g\rangle_0  = $ 
$ \langle f,\hH^*(g)\rangle_0$. Here $f,g$ are either
from $\v_0$ or from $\v_1$. All products of $\{X_b,Y_b, T_j,
\pi_r, q,t_\nu\}$ are unitary operators.

b) The polynomials $\{e_b,b\in B\}$ are eigenvectors of
 the operators $\{L_f\equal f(Y_1,\cdots, Y_n), f\in \C[x]\}$:
$$
\eqalignno{
&L_{\bar{f}}(e_b)\ =\ f(\#b)e_b, \where
\#b\equal\pi_b=b\om_b^{-1},
&\eqnu
\label\Yone\eqnum*
\cr
& x_a(\#b)\equal x_a(q^b t^{-\om_b^{-1}(\rho)})\ =\ 
q^{(a,b)}\prod_\nu t_\nu^{-(\om_b^{-1}(\rho_\nu),a)},\ w\in W.
&\eqnu
\label\xaonb\eqnum*
}
$$
\endproclaim
\label\YONE\theoremnum*
{\it Proof.} Assertion a) for $\v_0$ is from [C2].
Using (\ref\vovone) we come to $\v_1$ (a formal proof is
equally simple). Since
  operators $\{Y_b\}$ are unitary
relative to $\langle\ ,\ \rangle_0$
and leave all
$\Si(a),\Si_*(a)$ invariant
(Proposition \ref\TONX), their eigenvectors in $\C_{q,t}[x]$ are
exactly $\{e\}$. See [M3,C4].
\proofbox

The theorem  results immediately
in the orthogonality of $\{e_b\}$  for pairwise
distinct $b$. Macdonald also gives the formula for the squares
of $e_b$ (for $t_\nu=q^k,\ k\in \Z_+$)
and writes that he deduced it from the corresponding 
formula in the  $W$-symmetric case (proved in [C2]). The general case
was considered in [C4] where we used the recurrence relations. 
 A direct simple proof (based on the intertwiners)
will be given below.

The 
{\it symmetric Macdonald polynomials} form  a basis
 in the space 
$\C_{q,t}[x]^W$ of all $W$-invariant polynomials 
and
can be expressed as follows:
$$
\eqalign{
&p_{b}\ =\ \p_b^t e_{b}\ =\ \p_b^1 e_{b},\ \ b=b_+\in B_+,\cr 
&\p^t\equal\sum_{c\in W(b)}
\prod_\nu t_\nu^{l_\nu(w_c)/2} \hT_{w_c},\
w_c\equal \om_c^{-1}w_0,\ \p_b^1=\p_b^{t=1}.
}
\eqnu
\label\symmetr\eqnum*
$$

This presentation is from [M3,C4] (from [O2] in the differential
case). Here one can take the complete symmetrizations (with proper
coefficients) since $e_b$ is $W_b$-invariant for the 
stabilizer $W_b$ of $b$.
Macdonald  introduced these polynomials in [M1,M2]
by the  conditions 
$$
\eqalignno{
&p_{b}-m_{b}\ \in\ \Si_+(b),\ \langle p_{b},\ m_{c}\rangle_0 = 0, 
\   b\in B_+,\ c\succ b,
&\eqnu
\label\macdsym\eqnum*
}
$$
for the monomial symmetric functions $m_b=\sum_{c\in W(b)}x_c$.
 One can 
also define $\{p\}$ as eigenvectors for the
 ($W$-invariant) operators $L_f,\ $  $f\in \C_{q,t}[x]^W$:
$$
\eqalignno{
&L_{f}(p_b)\ =\ f(q^{b^o}t^{\rho} ) p_b,\
b\in B_+,b^o=-w_0(b),
&\eqnu
\label\Lf\eqnum*
}
$$
normalized as above.
Applying  any elements from $\h_Y=<T_j,Y_b>$
to   $e_c\ (c \in W(b_+))$
we get solutions
of (\ref\Lf), because
 symmetric $Y$-polynomials are central
in $\h_Y$ (due to I. Bernstein).
It readily gives the coincidence of (\ref\symmetr) and
(\ref\Lf).

%\vfill%
{\bf Functional representations}. The representations 
$\f_\xi, \De_\xi$ also have invariant skew-symmetric
forms.  Let
$$
\eqalignno{
&\mu_1(bw) = \mu(bw)/\mu(1)\equal 
&\eqnu
\label\muone\eqnum*
\cr
 \prod_{a \in R_+^\vee}
\prod_{i=0}^\infty &{(1-x_a(bw)q_a^{i}) (1-x_a^{-1}(bw)q_a^{i+1})
(1-x_a(1) t_aq_a^{i}) (1-x_a^{-1}(1)t_a q_a^{i+1})
\over
(1-x_a(1)q_a^{i}) (1-x_a^{-1}(1)q_a^{i+1})
(1-x_a(bw) t_aq_a^{i}) (1-x_a^{-1}(bw)t_a q_a^{i+1})}.
}
$$
Here $bw\in W^b \ni id=1,\ q_a=q_{\nu}=q^{2/\nu} \for \nu=\nu_a,
\ x_a(bw)=x_a(q^b w(\xi))$
(see (\ref\deltag)).

We ignore the convergence problem because  $\mu_1(bw)\in \C_\xi$
(actually it belongs to $\Q(q_\nu,t_\nu, \xi_i^{2/\nu_i})$: 

\proclaim{Proposition}  
a) Using $ \la(\hw) \for \hw=bw\in W^b$ from (\ref\lambda),
$$
\eqalignno{
& \mu_1(\hw) =\mu_1(\hw)^* \ = \
\prod_{[\al,j]\in \la(\hw)}
\Bigl(
{
t_\al^{-1/2}-q_\al^jt_\al^{1/2} x_a(\xi)\over
t_\al^{1/2}-q_\al^jt_\al^{-1/2} x_a(\xi)
}
\Bigr),
&\eqnu
\label\muval\eqnum*
}
$$
where $a=\al^\vee$, and 
we extend the conjugation $^*$ from $\C_{q,t}$
to $\C_\xi$ setting
$ \xi_i^*\ =\  \xi_i^{-1}$.

b) The following $\C_\xi$-valued
scalar product is well-defined for $f,g$
from the $\HH$-submodule of finitely supported
functions  $F_\xi\subset \f_\xi=\hbox{Funct}(W^b,\C_\xi)$:
$$
\eqalignno{
&\langle f,g\rangle_1\ = \sum_{\hw\in W^b}
\mu_1(\hw) f(\hw)\ g(\hw)^*\ =\ 
\langle g,f\rangle_1^*.
&\eqnu
\label\innerdel\eqnum*  
}
$$

c) Assertion a) from Proposition \ref\YONE holds
for $F_\xi$ and  $\De_\xi$, where the latter module
is endowed with the scalar product 
$$
\eqalignno{
&\langle f,g\rangle_{-1} = \sum_{\hw\in W^b}
(\mu_1(\hw))^{-1} u_{\hw} v_{\hw}^*,\ 
f=\sum u_{\hw}\de_{\hw}, g=\sum v_{\hw}\de_{\hw}.
&\eqnu
\label\innerdmin\eqnum*  
}
$$
Namely, $\langle H(f),g\rangle_{\pm 1}  = $ 
$ \langle f, H^*(g)\rangle_{\pm 1}$.
\endproclaim

{\it Proof.} Since $x_{\ta}(\hw)=x_{\ta'}(1) \for
\ta=\tal^\vee\in (R_+^a)^\vee, \where 
\ta'\equal\hw^{-1}(\ta)$, one has for $\hw=bw$:
$$
\eqalign{
\mu_1(\hw)\ =\ &\prod_{\tal \in R_+^a}
 { (1-x_{\ta}(\hw))(1-t_a x_{\ta}(1))
\over
 (1-x_{\ta}(1))(1-t_a x_{\ta}(\hw)) }
 \cr
=\  &\prod_{\tal \in R_+^a}
 {(1-x_{\ta'}(1))(1-t_a x_{\ta}(1))
\over
 (1-x_{\ta}(1))(1-t_a x_{\ta'}(1)) } \cr
=\  &\prod_{\tal \in \la(\hw)}
 {(1-x^{-1}_{\ta}(1))(1-t_a x_{\ta}(1))
\over
 (1-x_{\ta}(1))(1-t_a x^{-1}_{\ta}(1)) }.
}
\eqnu
\label\muaff\eqnum*
$$
Here we use that $\hw^{-1}(R_+^a)=\{-\la(\hw)\}\cup
\{R_+^a\setminus \la(\hw)\}$.
The invariance of $\mu_1(\hw)\in \C_\xi$ with respect
to the conjugation $^*$ is obvious. Other
statements  are completely analogous to those
for $\mu_0$ (and follow from them).
The key relation 
$$
\hH\mu(X)\ =\ \mu(X)(\hH^{*})^{+}, \ H\in \HH\ ,
\eqnu
$$
\label\deHdem\eqnum*
readily holds after the discretization. Here 
 by $^+$ we mean the anti-involution  
$$\hw^+= \hw^{-1}\in W^b, \ x_b^+= x_b^{-1}, b\in B,
 \ q,t\to q^{-1},\ t^{-1}.
$$
Its discretization conjugates the values of functions
from $\f_\xi$ and the coefficients of $\de_{\hw}$ in
$\De_\xi$ (fixing $\de_{\hw}$).
\proofbox

The {\it characteristic functions}
$f_{\hw}\in F_\xi (\hw\in W^b)$ are defined from the 
relations $\ f_{\hw}(\hu)=$ $\de_{\hw,\hu}\ $ for the
Kronecker delta. The action of 
the operators
$X_b$ on them is the same as for $\{\de_{\hw}\}$:
$$ X_b(f_{\hw})=x_b(\hw)f_{\hw},\ 
   X_b(\de_{\hw})=x_b(\hw)\de_{\hw},\ b\in B, \hw\in W^b.    
$$
Moreover the map  
$$
\eqalignno{
&f_{\hw}\ \to\ \mu_1(\hw)\de_{\hw},\ \hw\in W^b,
&\eqnu
\label\isodechar\eqnum*
}
$$
establishes an $\HH$-isomorphism between  $F_\xi$ and  $\De_\xi$, 
taking $\lan\ ,\ \ran_1$ to $\lan\ ,\ \ran_{-1}$.  
It readily  results  from the formulas:
$$
\eqalignno{
{}^\de T_i(f_{\hw})\ =\ 
{t_i^{1/2}x_{a_i}^{-1}(w(\xi))q^{-(a_i,b)} - t^{-1/2}\over
{x_{a_i}^{-1}(w(\xi))q^{-(a_i,b)}- 1}  }\
&f_{s_i\hw}\cr
-{t_i^{1/2}-t_i^{-1/2}\over
{x_{a_i}(w(\xi))q^{(a_i,b)}- 1}  }\
&f_{\hw} \for 0\le i\le n,
&\eqnu
\label\Tfchar\eqnum*
}
$$
and the formulas for the action of $\{\pi_r\}$.

Let us consider the special case $\xi=t^{-\rho}$ (see (\ref\deltaf)). 
Using the pairing $(\ref\depair)$, we see that the subspace
$$
F^\#=\oplus_{\hw\not\in \#B}\C_{q,t}f_{\hw}\subset
F(-\rho)=F_{t^{-\rho}}, 
\eqnu
\label\forho\eqnum*
$$
where $\#B=\{\#b=\pi_b\in W^b,\ b\in B\}$,
is an $\HH$-submodule. It is exactly the radical of the
form $\lan\ ,\ \ran_1$, which is well-defined for such $\xi$.
Indeed, any $\hw$ can be uniquely represented
in the form (see [C2])
$$
\hw=\pi_b w, \where b=\hw\lan b\ran,\ w\in W,\
l(\hw)=l(\pi_b)+l(w).
$$ 
Hence, $\{\hw\not\in \#B\}\Rightarrow
\{\al_i\in \la(\hw) \hbox{\ for \ some\ } i>0\}
\Rightarrow \{\mu_1(\hw)=0\}$. On the other hand,
$$
\eqalignno{
&\mu_1(\#b)\ =\ 
\prod_{a\in R_+^\vee}
\Bigl(
{
t_\al^{-1/2}-q_\al^jt_\al^{1/2} x_a(t^\rho)\over
t_\al^{1/2}-q_\al^jt_\al^{-1/2} x_a(t^\rho)
}
\Bigr), 
&\eqnu
\label\muorho\eqnum*
}
$$
where the product is over the set $\la'(\pi_b)=
\{[\al,j], [-\al,j]\in \la(\pi_b)\}$. Explicitly (see (\ref\lambpi)),
$$
\eqalignno{
\la'(\pi_b)=\{[\al,j], 0<j<-(\al,b_-)\if &(\al,b)>0,\cr
 0<j\le -(\al,b_-)\if &(\al,b)<0\},\ b_-=\om_b(b).  
&\eqnu
\label\jbset\eqnum*
}
$$
Since  $j>0$ in either case, $\mu_1(\#b)\neq 0$.
The  map (\ref\isodechar) identifies 
 $F_\#\equal F(-\rho)/F^\#$ with $\De_\#$.

Generally speaking, the problem is to go from $F_\xi$ to
$\f_\xi$, for instance,
to introduce and decompose the module of all square
integrable functions. At least, one can try to figure out which
functions form the image of $\v_0$ in $\f(-\rho)/\f^\#$
are square integrable. We will touch upon this problem  in
the next section.

%
%
%		Section 5
%
%
%\vskip 10pt
\section { Applying intertwiners } 
Here we will use the intertwiners as creation operators for
the nonsymmetric polynomials and establish connections with
the represenation $\De_\#$.
We assume in this section that $\xi=t^{-\rho}$ and 
$x_a(\#b)=x_a(\pi_b)=x_a(q^b t^{-\om_b^{-1}(\rho)})$.
The notations are from the previous sections. Let us set:
$$
\eqalignno{
&\Phi_j^b\ =\ \Phi_j(\#b)\ =\
T_j + (t_j^{1/2}-t_j^{-1/2})(x_{a_j}(\#b)-1)^{-1},
&\eqnu
\label\Phijb\eqnum*
\cr
&G_j^b\ =\ (\Phi_j\phi_j^{-1})(\#b)\ =\
{ T_j + (t_j^{1/2}-t_j^{-1/2})(x_{a_j}(\#b)-1)^{-1}
\over
t_j^{1/2} + (t_j^{1/2} -t_j^{-1/2})(x_{a_j}(\#b)-1)^{-1} }.
&\eqnu
\label\Gjb\eqnum*
}
$$
\proclaim {Main Theorem}
Given $c\in B,\ 0\le j\le n$ such that $(\al_j, c+d)> 0,$ 
$$
\eqalign{ 
&e_{b}q^{-(b,b)/2} \ =\ t_j^{1/2}\Phi_j^c\llan e_c \rran q^{-(c,c)/2}  
\for b= s_j\lan c\ran, \cr
&\hde_{b}\ =\ G_j^c(\hde_c),\
\hf_{b}\ =\ \tG_j^c(\hf_c), \ \tG_j^c=(\phi_j^{-1}\Phi_j)(\#c), 
}
\eqnu
\label\Phieb\eqnum*
$$
where $\hf_b$ is the image of $f_{\#b}$ in $F_{\#}= F(-\rho)/F^{\#}$,
$\hde_b=\de_{\#b}$.
This inequality is equivalent to $(\al_j, b+d)< 0$ and to the relation
 $l(\pi_b)=l(\pi_c)+1$.   If $(\al_j, c+d)=0$ ( $\Leftrightarrow
\phi_j(c)=0$) then
$$
\eqalignno{
&T_j \llan e_c \rran \ =\ t_j^{1/2} e_c, \ 0\le j\le n, 
&\eqnu
\label\Tjeco\eqnum*
}
$$
which  gives that $s_j(e_c)=e_c$ when
$j>0$.  Also for any $b=\pi_r\lan c\ran,\ r\in O$,
$$
\eqalignno{
&q^{-(b,b)/2}e_b= q^{-(c,c)/2}\pi_r\llan e_c\rran,\ 
\hf_{b}=\pi_r(\hf_{c}),\  \hde_{b}=\pi_r(\hde_{c}).
&\eqnu
\label\pireb\eqnum*
}
$$
\endproclaim
{\it Proof.}
The element $\tau^{-1}\vep\tau(\Phi_j^c)=\eta(\Phi_j^c)=\Phi_j^c$ 
sends $e_c$ considered
as an element of $\v_1$ to a 
nonzero polynomial 
proportional to $e_b$. Here we use that $\Phi_j$ is an
$X$-intertwiner and $\phi_j(\#c)\neq 0$. The latter results 
from the inequality $(\al_j, c+d)> 0$ (see (\ref\aljb)).
The coincidence of $\eta(\Phi_j^c)$ and $\Phi_j^c$ is due to
(\ref\eta).  
To make this  more obvious one can involve
$\ga$. Then $\Phi_j^c\llan e_c\rran=
\ga \Phi_j^c(e_c\ga^{-1})$ and the reasoning  
gets rather straightforward.

Thus $\Phi_j^c\llan e_c\rran= u e_b q^{\{(c,c)-(b,b)\}/2}$.
The leading term (up to $x_{b'}, b'\in\si_*(b)$) 
of the second expression is $u x_b$. We need to find $u$.
Setting  $e_c = x_c+\sum_{a\in \si_*(c)} v_{ca} x_a$,
the elements $\Phi_j^c\llan x_a\rran$ do not contribute 
to $x_b$ if $(\al_j, a+d)\le 0$ (see (\ref\tlonx)).
We pick a minimal element $a'=s_j\lan a \ran$ 
realtive to the partial ordering $\succ$ from the set 
$$
S'=s_j\lan S \ran 
 \for
S= \{a,\ v_{ca}\neq 0,\ (\al_j, a+d)> 0\}.
$$ 
Due to (\ref\tmonx),  $\Phi_j^c\llan x_{a}\rran= v_a x_{a'}$ 
for $v_a\neq 0$ modulo $\Si_+(a')$. If $a'\neq b$ then 
 $\Phi_j^c\llan x_c\rran$ contains $x_{a'}$ with nonzero
coefficient for $a'\not\in \si(b)$ which is impossible.
Hence, $S'\subset \si(b)$ and  any elements
$c \neq a\in S$ go to $b\neq a'\in \si(b)$. It gives that
 $u=t_j^{-1/2}$. 

Similarly, 
$ T_j \llan e_c \rran \ =\ u_j e_c \ (0\le j\le n)\ $ 
if $(\al_j,c+d)=0$, and we can apply the same argument to
see that $u_j=t_j^{1/2}$.
The statements about the characteristic and delta-functions
are checked by simple direct calculations.
\proofbox

We can reformulate the Theorem in the following way.
Let us introduce the {\it renormalized polynomials}:
$$
\eqalignno{ 
&\he_{b}\ =\ (\pi_r G_l^{c_l}\ldots G_1^{c_1})\llan 1 \rran,
\where  
&\eqnu
\label\ehatb\eqnum*
 \cr
&c_1=0, c_2=s_{i_1}\lan c_1\ran,\ldots,
c_l=s_{i_l}\lan c_{l-1}\ran, \for \pi_b= \pi_r s_{j_l}\ldots s_{j_1}.
}
$$
They are well-defined,
do not depend on the particular choice of
the decomposition of $\pi_b$ (not necessarily
reduced), and are proportional to $e_b$ for all
$b\in B$. The coefficients of proportionality
 (always nonzero) can be readily calculated using (\ref\Phieb):
$$
\eqalignno{ 
e_{b}q^{-(b,b)/2} \ =&\ \prod_{1\le r\le l}
\bigl(t_{j_r}^{1/2}\phi_{j_r}(c_r)\bigr)\ \he_b\cr
=&\ \prod_{[\al,j]\in \la'(\pi_b)}
\Bigl(
{ 
1- q_\al^{j}t_\al x_a(t^\rho)
 \over
1- q_\al^{j}x_a(t^\rho)
}
\Bigr)\ \he_b,
&\eqnu
\label\ebebhat\eqnum*
}
$$
where $\la'(\pi_b)$ is from (\ref\jbset).

\proclaim {Corollary}
a) Let $\chi: \v_1\to \De_{\#} $ be a $\ast$-linear isomorphism
sending $\sum u_b\he_b\to \sum u_b^*\hde_{b}$ for $u_b\in \C_{q,t}$. 
Then (see (\ref\eta)) 
$$
\eqalignno{ 
&\chi(H\llan v\rran)\ =\ \de(\eta(H))(\chi(v)) \for  v\in \v_1, H\in \HH.
&\eqnu
\label\isovd\eqnum*
}
$$

b)
The $\ast$-linear 
isomorphism $\be: \v_0\to \De_{\#}$ 
sending 
$$
\eqalignno{ 
&\sum u_b\he_b\to 
\sum u_b^{*}\ga(\#b)\hde_{b},\
\ga(\#b)=q^{(b,b)/2}x_b(t^{-\om_b^{-1}(\rho)})
&\eqnu
\label\isofour\eqnum*
}
$$
satisfies the relation:
$$
\eqalignno{ 
&\be(H(v))\ =\ \de(\vep(H))(\be(v)) \for  v\in \v_0, H\in \HH.
&\eqnu
\label\fourpol\eqnum*
}
$$
\endproclaim
\label\ISOVD\theoremnum*
{\it Proof.} This map satisfies (\ref\isovd) for 
$H\in \h_Y=<T_i,\pi_r>$ due to the previous considerations.
On the other hand,
$X_i(\hde_b)\ =\ x_i(\#b)\hde_b$ and 
$$
\tau^{-1}(Y_i)(\he_b)\ =\  (\ga^{-1}(Y_i)\ga)(\ga^{-1}e_b)
\ =\ x_i^{-1}(\#b)\he_b.
$$  
Since $\eta(X_i)=\tau^{-1}(Y_i)$ and $\vep: q\to q^{-1},t\to t^{-1}$,
the relation holds for all $\tau^{-1}(Y_i)$.
Hence it is always true (the $\eta$-images of $\tau^{-1}(Y_i),T_j,\pi_r$
generate the whole \HH\ ).

The second statement follows from the first, since $\vep=\tau\eta\tau^{-1}$
and $\tau$ acts as conjugation by the Gaussian $\ga$. The multiplication
by $^{\de}\ga$ on the delta-functions leads exactly
to the coefficients from (\ref\isofour) (see (\ref\degauss)).
\proofbox
 
Assertion b) is nothing else but the calculation of the 
Fourier transform of the nonsymmetric polynomials
up to a common factor. The  Fourier
transform from [C3,C4] induces $\vep$ on the operators (it
is a defining property).
Hence it is proportional to $\be$
(cf. [C4] using Corollary \ref\EVE below).
This construction also 
establishes a connection of the representations of the
affine Hecke algebra in $\C[X]$ and that in functions
on $I\backslash G/K$ (see the end of Section 3).
Both representations can be defined $p$-adically
and look very similar. However
it seems that a natural connection of these two representations
exists only at 
level of the $q,t$-theory (the question was
suggested to the author by D. Kazhdan).

\proclaim {Corollary}
The coefficients of the polynomials 
$$
\eqalignno{ 
 \prod_{[\al,j]\in \la'(\pi_b)}
\bigl(
1- q_\al^{j}x_a(t^\rho)\bigr)\ e_b, \
&\prod_{[\al,j]\in \la'(\pi_{b_+})}
\bigl(
1- q_\al^{j}x_a(t^\rho)\bigr)\ p_{b_+}, 
&\eqnu
\label\integ\eqnum*
\cr
 q^{(b,b)/2}&\prod_{[\al,j]\in \la'(\pi_b)}
\bigl(
1- q_\al^{j}t_\al x_a(t^\rho)\bigr)\ \he_b 
&\eqnu
\label\integp\eqnum*
}
$$
belong to $\Q[q^{\pm 1},t_\nu^{\pm 1}]$.
\endproclaim
\label\INTEG\theoremnum*
{\it Proof.} We use (\ref\symmetr) for the symmetric
polynomials $p_{b_+}$.
\proofbox

\proclaim {Corollary}
For $b\in B, \ b_-=\om_b(b)\in B_-$, and $^{\de}\ga$
from (\ref\degauss),
$$
\eqalignno{ 
& \he_b(\#)=\he_b(t^{-\rho})\ =\ 
^{\de}\ga^{-1}(\pi_b)\ =\ t^{(\rho,b_-)}q^{-(b,b)/2}.
&\eqnu
\label\eve\eqnum*
}
$$
\endproclaim
\label\EVE\theoremnum*
{\it Proof.}
It follows immediately from the Main Theorem of [C4]
and (\ref\ebebhat).
We will  outline a direct reasonning which also gives
another proof of the evaluation formula for nonsymmetric polynomials.

Let us start with the following symmetric, non-degenerate
pairing from [C4]
on $f,g\in \C_{q,t}[x]$:
$$
\eqalignno{
 &[\![f,g]\!]_0\ =\  \{L_{\bar{f}}(g(x))\}({\#}),\
&\eqnu
\label\Fourier\eqnum*
\cr
&\bar{x}_b\ =\ x_{-b}\ =\ x_b^{-1},\ \bar{ q}\ =\  q,\ 
\bar{t}\ =\ t, 
}
$$
where
$L_f$ is from Proposition \ref\YONE.
The corresponding anti-involution is the composition $\vph= 
\vep\ast=\ast\vep$, sending $X_i\to Y_i^{-1}$ and preserving
$T_i, q,t\ (1\le i\le n)$:
$$
\eqalign{
 &[\![H(f),g]\!]_0\ =\ [\![f,{\vph(H)}(g)]\!]_0,\ H\in \HH.
}
\eqnu
\label\phiinvar\eqnum*
$$
Similarly, we can define the symmetric pairing 
$$
\eqalignno{
 &[\![f,g]\!]_1\ =\  \{(\ga^{-1}L_{\bar{f}}\ga)(g(x))\}(\#),\where
&\eqnu
\label\Fourone\eqnum*
\cr
&[\![H(f),g]\!]_1\ =\ [\![f,{\psi(H)}(g)]\!]_0 \for
\psi=\tau^{-1}\vph\tau=\ast\eta=\eta\ast,
}
$$
and the same action $H(f)=\hH(f)$.
The claim is  that for the standard pairing $\{\ ,\ \}$
from (\ref\depair),
$$
\eqalign{
 &[\![f,g]\!]_1\ =\ \{f,\ze(g)\},
}
\eqnu
\label\fgchi\eqnum*
$$
where $\ze$ is a $\C_{q,t}$-linear counterpart of $\chi$:
$$
\eqalignno{
&\ze(\sum_b u_b\he_b)\ =\  \sum_b u_b \hde_{b} \for u_b\in \C_{q,t}.
&\eqnu
\label\zeiso\eqnum*
} 
$$
Indeed, the anti-involutions corresponding to both pairings 
fix all $T_j,\pi_r$. 

Let us rewrite (\ref\Fourone) as follows:
$$
\eqalignno{
 &[\![f,g]\!]_1\ =\  \{(L_{\bar{f}}(\ga g(x))\}({\#})\ =\ 
\{(L_{\vph(\ga)(\bar{f})}( g(x))\}({\#}).
&\eqnu
\label\Fourga\eqnum*
}
$$
Recall that $\ga$ is normalized by the condition $\ga(\#)=1$.
We do not give a formal definition of $\vph(\ga)$.
Anyway $\vph(\ga)(e_b)=\ga(\#b)(e_b)$ as it is for polynomials
of $Y$. Therefore
$[\![\he_b,\he_c]\!]_1 = \ga(\#c)\he_b(\#c)\he_c(\#)$.
Since it also coincides with $\{\he_b,\hde_c\}=\he_b(\#c)$, 
we come to the required formula.
\proofbox

Using the same argument, we can get interesting  recurrence formulas
for the values of $\he_b$. A typical example is the following

\proclaim{Corollary}
Let as assume that $(\al_j, c+d)=-(\al_j,b+d)$ for $0\le j\le n,\
c,b\in B$. Then 
$$
\eqalignno{
 &\he_b(\# c)\ =\ \{\he_b,\hde_c\}\ =\
\{\he_{b'},\hde_{c'}\}\ =\ \he_{b'}(\# c'),
&\eqnu
\label\valrecur\eqnum*
}
$$
where $b'=s_j\lan b\ran,\ c'=s_j\lan c\ran$.
\endproclaim
\proofbox

Following the same lines, we will calculate the norms of
$\he_b$. Actually we can use again the Main Theorem from
[C4], since the coefficient of proportionality
 $e_b/\he_b$
has been already known. However the construction below
establishes an important direct connection with $\mu_1$. This can
be deduced from [C4] as well 
but in a more complicated way.

\proclaim{ Theorem} For $b,c\in B$ and the Kronecker delta $\de_{bc}$,
$$
\eqalignno{
 &\lan \he_b,\he_c\ran_0\ =\ \lan \de_b,\de_c\ran_{-1}\ =\ 
\de_{bc}\mu_1^{-1}(\#b)\cr 
&=\ \de_{bc}\prod_{[\al,j]\in \la'(\pi_b)}
\Bigl(
{
t_\al^{1/2}-q_\al^jt_\al^{-1/2} x_a(t^\rho)\over
t_\al^{-1/2}-q_\al^jt_\al^{1/2} x_a(t^\rho)
}
\Bigr). 
&\eqnu
\label\normhats\eqnum*
}
$$
\label\NORMSHAT\theoremnum*
\endproclaim
{\it Proof.} 
We claim that 
$$
\eqalign{
&\lan f,g\ran_0\ =\ \lan \ze(f),\ze(g)\ran_{-1},\where
f,g\in \C_{q,t}[x],
}
\eqnu
\label\zenorms\eqnum*
$$
for the map 
$\ze$ from (\ref\zeiso)
(or $\chi$ from Corollary \ref\ISOVD, doesn't matter).
We represent
$\pi_b\ =\ \pi_r s_{j_l}\ldots s_{j_1}$ (see (\ref\ehatb)), 
$$
\eqalignno{ 
&\he_{b} = (\pi_r G_l^{c_l}\ldots G_1^{c_1})\llan 1 \rran,\
\hde_{b} = (\pi_r G_l^{c_l}\ldots G_1^{c_1})(\hde_0).
}
$$
Then we use the relations
$G_j^*=G_j,\ \pi_r^*=\pi_r^{-1},\ \ga^*=\ga^{-1}$:
$$
\eqalignno{ 
&\lan\he_{b},\he_c\ran_0\ =\ \lan\he_{b}\ga^{-1},\he_c\ga^{-1}\ran_0\ =\cr
&\lan\ga^{-1},  G_1^{c_1}\ldots G_l^{c_l}\pi_r^{-1}(\he_c\ga^{-1})\ran_0\ =
\ \lan\ 1,  G_1^{c_1}\ldots G_l^{c_l}\pi_r^{-1}
\llan \he_c\rran\ \ran_0,
&\eqnu
\label\squaro\eqnum*
\cr
&\lan\de_{b},\de_c\ran_{-1} =
\lan \de_0, G_1^{c_1}\ldots G_l^{c_l}\pi_r^{-1}(\de_c)\ran_{-1}.
&\eqnu
\label\squarde\eqnum*
}
$$
Now (\ref\zenorms) is obvious since
 $\ze$ leaves $\{G_j^c$, $\pi_r\}$ invariant and
$\lan 1,1\ran_0=1=\lan \hde_0,\hde_0\ran_{-1}$.
\proofbox 
 
Here we can replace $\lan\ ,\ \ran_0$ by any scalar
product on polynomials providing the $\ast$-invariance and
the normalization $\lan 1,1\ran_0=1$. 
Indeed, 
$ G_1^{c_1}\ldots G_l^{c_l}\pi_r^{-1}\llan \he_c\rran $
 can be linearly expressed in terms of proper $\he_a$.
Because different $\he_b$ are pairwise orthogonal, 
 (\ref\squaro) equals the coefficient of $\he_0=1$
in this expression.

\proclaim{ Proposition}
a) Given $c\in B$, let us assume that the infinite sums for the
scalar products $\lan ^\de\he_a, ^\de\he_b\ran_1$ are
absolutely convergent for any $B\ni a,b \succeq c$.
Then
$$
\eqalignno{ 
&\lan ^\de\he_a, ^\de\he_b\ran_1\ =\ 
\lan 1, 1\ran_1 \lan \he_a, \he_b\ran_0, \ B\ni a,b \succeq c.
&\eqnu
\label\jacknorm\eqnum*
}
$$

b) We suppose that $|q|\neq 1, \
q_\al^j x_{\al^\vee}(\xi)\neq t_\al$
for all $[\al,j]\in R_+^a$ (cf. (\ref\indsph)), 
and $c=c_-\in B_-$. 
The necessary and sufficient 
condition for the absolute convergence of  
$\lan ^\de\he_a, ^\de\he_b\ran_1$ for all $a,b$ above 
is as follows.
Setting 
$$t_i=q_i^{k_i},\ c_+=w_0(c_-)\in B_+,\  
r_k=\sum_{i=1}^n k_i b_i, \ \Re(r_k)=\sum_{i=1}^n \Re(k_i) b_i,
$$
where $\Re$ is the real part,
$$
\eqalignno{ 
& 2\Re(r_k)+c_+ -c_-\ =\ \sum_{i=1}^n p_i a_i \hbox{\ s.t.\ } \R\ni p_i<0.
&\eqnu
\label\jackcond\eqnum*
}
$$
\label\JACKNORM\theoremnum*
\endproclaim

{\it Proof}. The first part results from the above remark.
The convergence  is checked  similar to [Ito] (use (\ref\muval)).
Given an  arbitrary $c\in B$,
the products $\lan ^\de\he_a, ^\de\he_b\ran_1$
are absolutely convergent if and only if
$$
\eqalignno{ 
& 2\Re(r_k)+w(a-b)\ =\ \sum_{i=1}^n p_i(a,b,w) a_i,\where  p_i(a,b,w)<0.
&\eqnu
\label\jackgcon\eqnum*
}
$$
for all $B\ni a,b \succeq c,\ w\in W$. The set $\{w(a-b)\}$
is especially simple when $c=c_-$. In this case $w(a-b)\le
c_+-c_-$. Since the latter difference also
belongs to this set, we come
to (\ref\jackcond).
\proofbox

The value of the constant $\lan 1, 1\ran_1$
is directly related to the Aomoto conjecture (see ibid.)
recently proved by Macdonald [M4].
It can be also calculated by
means of the discretization of the shift 
operators from [C2]. One arrives at simple relations
connecting $\lan 1,1\ran_1$ for $k$ and $k-1$ and then can proceed as
in [O1]
(the differential case).  

Replacing $\mu$ by its $W$-symmetric counterpart (due to Macdonald)
$$
{\mu'}\ =\ \mu\prod_{a\in R_+^a}(1-x_a^{-1})(1-t_a x_a^{-1})^{-1},
\eqnu
\label\symmu\eqnum*
$$
we can introduce ${\mu'}_1={\mu'}/{\mu'}(1)$ following
(\ref\muone) and 
${\mu'}_0={\mu'}/\lan {\mu'}\ran_0$. 
Then providing the conditions
from b),
$$
\eqalignno{ 
&\lan p_{a_+},p_{b_+}\ran'_1 = 
A_\xi \lan p_{a_+},p_{b_+}\ran'_0,\
A_\xi=\sum_{d\in B}{\mu'}_1(q^d\xi),
&\eqnu
\label\jacksym\eqnum*
\cr 
&\lan f, g\ran'_1 \ =\  
\sum_{d\in B} {\mu'}_1(q^d\xi)
f(q^d\xi)g(q^d\xi),\cr 
&\lan f, g\ran'_0\ =\ \lan {\mu'}_0(x) f(x)g(x^{-1})\ran_0.
}
$$
The coefficient of proportionality $A_\xi$  is right from the
Aomoto conjecture. 
The proof is based on (\ref\jacknorm) and Proposition 4.2
from [C2]. Both forms make the $L$-operators 
(see (\ref\Lf)) self-adjoint.
These pairings work well for symmetric polynomials
$f,g$ only but have some merits because the summation is over $B$ and
it is not necessary
to conjugate $q,t,\xi$. Here we also should
assume that $f,g$ are ``real''
with respect to this conjugation ($p_{b_+}$ are real).

We note that our approach generalizes the
calculation of the norms of Opdam's nonsymmetric polynomials
from the same $W$-orbit  [O2]. He  used
the intertwiners too but only non-affine ones (in the differential
case). There are quite a few papers on Jackson integrals
of the $q$-polynomials. Mostly they are one-dimensional.
Let us mention a recent work [StK].

%
%
%		Section 5
%
%
%\vskip 10pt
\section { Induced and co-spherical representations} 

We   return to the  case of general $\xi$. In
this section we  treat  $q, t_\nu, \xi_i$ as nonzero complex numbers.   
The delta-functions
$\de_{\hw}$ will be considered 
as  characters of $\C[x]$ using the pairing
(\ref\depair):
$$
\eqalignno{ 
& \de_{bw}(x_{\ta})\equal x_{\ta}(bw)\ =\ q^{k+(a,b)}x_{w^{-1}(a)}(\xi),
\where \ta=[a,k]. 
&\eqnu
\label\dechar\eqnum*
}
$$ 
They form the $W^b$-orbit of $\de_1=\xi$. Different
$\de_{\hw}$ can coincide
for certain $q,\xi$. We will identify them in this case.

From now on $q$ is not a root of unity. 
This hypothesis is necessary and sufficient to make 
the  {\it stabilizers} $W_{\xi}^b(\hu)\equal\{\hw\in W^b,
\de_{\hw\hu}=\de_{\hu}\}$ finite for all $\hu$.
We will also assume that there exists a {\it primitive}
character $\xi_o=\de_{\hu_o}$ such that $W^o\equal W_{\xi}^a(\hu_o)$ 
is generated by the elements from
$S^o\equal W^o\cap\{s_0,\cdots,s_n\}$. It means that
$W^o$ is the Weyl group of the non-affine Dynkin graph $\Ga^o\in\Ga^a$
(not necessarily connected)  with $S^o$ as the set of vertices.
The existence of the primitive character in the $W^b$-orbit of $\xi$
 always holds true for rather general $q,t$
(say, for degenerate double Hecke algebras below).
We denote the non-affine Hecke algebra corresponding to $\Ga^o$
by $\H^o$. Adding $\{ X_{a_i}, \ s_i\in S^o\}$ to $\H^o$ one gets the
affine Hecke algebra $\h_X^o$, which will be 
considered as a subalgebra
of $\HH$. 

The functional representations $\f_\xi$ and $\De_\xi$
can be introduced using the same formulas if 
$$
\eqalignno{ 
&
x_{\al^\vee}(\xi)q_\al^j \neq 1 \for
\hbox{all\ } \al\in R, j\in \Z. 
&\eqnu
\label\neqone\eqnum*
}
$$
All $\de_{\hw}$ are linearly independent ($W^o=\{1\}$)
and the pairing with $\v_0$ remains nondegenerate
in this case.
Functional representations can be defined without imposing 
(\ref\neqone), but we will not discuss it here.

An $\HH$-module $I$ is {\it co-spherical} 
if it  contains no
 submodules $V\neq \{0\}$ such that all $\h_Y$-invariant homomorphisms
$\om: I\to \C$ vanish on $V$.
By the $\h_Y$-invariance, we mean that 
$$
\eqalignno{ 
&\om(T_j (v))=t_j^{1/2}v,\
\om(\pi_r(v))=v \for 0\le j\le n,\ r\in O,\ v\in I.
&\eqnu
\label\cospher\eqnum*
}
$$
We will consider only  modules where the set  $\{\om\}$ 
contains not more than one element  adding this to the 
definition of co-spherical modules. Generalizing,
the invariance condition (\ref\cospher) can
be introduced for any character of $\h_Y$. The corresponding
extension of the results below  is straightforward.

The module  $\De_\xi$ is  co-spherical.
Indeed, $\De_\xi\ni f$ posses only one invariant homomorphism
$\om(f)=\{1,f\}$. There are no  $\HH$-submodules $V$  
 such that 
$\{1, V\}=0$, because $\v_0$ is $\HH$-generated by $1$
and the pairing $\{\ ,\ \}$ is nondegenerate. 
We also note that $\De_\xi$ is isomorphic to $\De$ defined
for any characters $\de_{\hw}$ taken instead of $\xi$.
The corresponding map is the right multiplication
by $\hw$ ($\de_{\hu}\to \de_{\hu \hw}$). 

Let us figure out when 
the {\it $X$-induced representations}  
are irreducible and co-spherical.
The definition is  standard (see e.g. [KL]). 
The induced representation
$I_\xi$ is the universal \HH-module 
generated by the element $v_\xi$ such that $X_i(v_\xi)=\xi_i v_\xi$.
As an $\h_Y$-module, it is isomorphic  to $\h_Y$ 
with the left regular action ($v_\xi$ is identified with $1\in \h_Y$). 
Hence, there exists only one 
invariant  $\h_Y$-homomorphism $\om$. It sends
 $T_j\to t_j^{1/2}, \pi_r\to 1$ after the identification.

\proclaim {Theorem}
a) The module $I_\xi$ is irreducible if and only if 
$$
\eqalignno{ 
&x_{\ta}(\xi)\neq t_{\al}^{\pm 1} \for \hbox{all\ }
\tal=[\al,k]\in R_+^a, \ \ta=\tal^\vee.
&\eqnu
\label\indirrep\eqnum*
}
$$
Irreducible $I_\xi$ are isomorphic for $\xi$ from the same
$W^b$-orbit (and only for such characters).

b) The module $I_\xi$ is co-spherical if and only if
$$
\eqalignno{ 
&x_{\ta}(\xi)\neq t_{\al}^{-1} \for \hbox{all\ }
\tal\in R_+^a.
&\eqnu
\label\indcosph\eqnum*
}
$$ 
If (\ref\indcosph) holds then induced modules
associated with  characters $\de_{\hw}$ satisfying the same
inequalities are isomorphic to $I_\xi$.

c) Under the same condition, $I_\xi$
 contains a unique nonzero irreducible submodule $U_\xi$.
It is co-spherical (i.e. has a nonzero $\om$).
Any other irreducible constituents are not co-spherical.
An arbitrary irreducible
co-spherical module possessing an eigenvector
with the character  $\de_{\hu}$ 
(i.e. belonging to the $W^b$-orbit of the
$\xi$ above) is isomorphic to $U_\xi$. 
\label\INDCOSPH\theoremnum*
\endproclaim

{\it Proof.}
The statements are parallel to the corresponding
affine ones. They result from the irreducibilty of
the simplest induced modules with the characters $\{1,\ldots,1\}$
by means of the technique
of intertwiners. See [KL] and the papers [Ro], [C5,C8].
We will mainly follow [C5]. 

Let us first renormalize $\Phi$ from (\ref\Phi),(\ref\Phijb) to avoid the
denominators:
$$
\eqalignno{
&\tPhi_j\ =\ 
T_j(X_{a_j}-1) + (t_j^{1/2}-t_j^{-1/2}),\ 0\le j\le n,
\cr
&\tPhi_j^{\hu} =\ \tPhi_j(\hu)\ =\
T_j(x_{a_j}(\hu)-1) + (t_j^{1/2}-t_j^{-1/2}),\ \hu\in W^b.
&\eqnu
\label\Phiind\eqnum*
}
$$
The corresponding $\tPhi_{\hw}$ are well-defined and 
enjoy the main property of the intertwiners (\ref\Phix).
The multiplication on the right
by the element $\tPhi_{\hu}(\hw)$
(which belongs to $\h_Y$)
is an $\HH$-homomorphism from the module $I_{\hw}\simeq \h_Y$
into the module $I_{\hw\hu}$.  Here
$I_{\hw}$
is the induced representation corresponding to the 
character $x_a\to x_a(\hw^{-1})\ =\  \de_{\hw^{-1}}(x_a)$.
Similarly, if $v$ is an $X$-eigenvector
corresponding to $\de_{\hu}$ then $\tPhi_{\hw}^{\hu}(v)$
is that associated with $\de_{\hw\hu}$.

We will use some general facts about eigenvectors and
induced representations (see [C5], Proposition 2.8, 
Lemma 2.10). They are based on the following definition
of the spaces of {\it generalized eigenvectors}:
$$
\eqalignno{
& I_\xi^N(\de_{\hu})\equal \{v\in I_\xi,\ 
(X_i-x_i(\hu))^N(v)\ =\ 0\for N\in \N\},\cr
&I_\xi^{\infty}(\de_{\hu})\ =\ I_\xi^N(\hu) \hbox{\ as\ } N\to \infty,\
I_\xi^{N}(\hu)\equal I_\xi^N(\de_{\hu}).
&\eqnu
\label\xeigen\eqnum*
}
$$
The spaces $\{I_\xi^{\infty}(\hu)\}$ are  finite dimensional and
$I_\xi=\oplus I_\xi^{\infty}(\hu)$ for pairwise different characters
$\de_{\hu}$. 

More generally,  $\HH$-modules $I$
such that  $I^{\infty}(\mu)$ are finite dimensional and 
$I=\oplus I^{\infty}(\mu)$ for all characters $\mu$
  of $\C[x]$  
constitute the {\it category $\o$} (see [BGG]). 
All irreducible modules
possessing $X$-eigenvectors belong to this category.  
If $I^\infty(\mu)\neq \{0\}$
then $I^1(\mu)\neq \{0\}$. Each $I\neq\{0\}$ contains at least one
nonzero irreducible submodule.

The first application is that
the set of all 
eigenvalues (i.e. the characters
associated with  $X$-eigenvectors) of $I_\xi$ is exactly
$\{\de_{\hw}\}$ for any $\xi$. Indeed, if sufficiently general
$\xi'$ tends to
$\xi$ then $I_{\xi}^\infty(\de_{\hu})$ is exactly the image of
the direct sum of $I_{\xi'}^\infty(\de_{\hu'})$ over all $\de_{\hu'}
\to \de_{\hu}$.

If $\hu_o$ is  primitive, then  
$$
\eqalignno{
& I_\xi^{\infty}(\hu_o)\ =\ \H^o\subset \h_Y.
&\eqnu
\label\infprim\eqnum*
}
$$
This space is an induced $\h_X^o$-module, corresponding
to the trivial character $\{X_{a_i}\to 1,\ s_i\in S^o\}$.
Moreover it is irreducible and contains a unique $X$-eigenvector
that is $1\in \H^o$. The proof of the last statement requires
some technique  (see Lemma 2.12 from [C5]). Using this, we see that 
the dimension of the
space $I^{1}(\hu_o)$ for
primitive $\hu_o$ is not more than one for
 any  submodules and subquotients ({\it constituents})
$I$ of $I_\xi$.

%\vskip 0.2cm
Let us check a). If $x_{\ta}(\xi)= t_{\al}^{\pm 1}$ for some $\tal$, then
there exists $\hw\in W^b$ such that the operator 
of the right multiplication by $\Phi_{\hw}(1)$ has a non-trivial
kernel. Hence $I_1=I_\xi$ is reducible. Let us assume that  
inequalities (\ref\indirrep) hold true.

All the elements $\tPhi_j^{\hu}\in \h_Y$ are invertible. 
Applying them to $v_\xi$ 
we can get an eigenvector $v_o$ corresponding to
$\de_{\hu_o}$. The latter generates the whole $I_\xi$ since
we can go back to $v_\xi$ using  $(\tPhi_j^{\hu})^{-1}$ in the
opposite order. Hence we can assume that $\hu_o=1$.

If $I_\xi$ contains an $\HH$-submodule $V\neq \{0\}$, then
there exists at least one $X$-eigenvector $v\in V$
(any vectors from $I_\xi$ belong to  finite dimensional
$X$-invariant subspaces). A proper chain of the intertwiners applied
to $v$ will produce a nonzero eigenvector $v_\xi'$ corresponding to
$\de_1=\xi$. It is proportional to $v_\xi$, since it
has to belong to $\H^o$ (see above). Hence $v$ generates $I_\xi$,
and $V=I_\xi$.

Let us come to b). 
If $ x_{\ta}(\xi)=t_{\al}^{-1}$  for a certain
$\tal\in R_+^a$, then using  a chain of invertible intertwiners
we can replace $v_\xi$ by an eigenvector $v\in I_\xi$ associated
to $\de_{\hw}$ such that $x_j(\hw)=\de_{\hw}(x_j)= t_{j}^{-1}$
for some index  $0\le j\le n$. Therefore we can assume, that  
$\xi_j=t_{j}^{-1}$. Then $\om=0$ on the nonzero $\HH$-submodule
$$
I_\xi\tPhi_j^1\ =\ (t_j^{1/2}-t_j^{-1/2})\h_Y(1-t_j^{-1/2}T_j).
$$

From now on,   $\xi$ will satisfy  (\ref\indcosph).
Applying the chains of intertwiners corresponding to
reduced decompositions of  elements $\hw$
to $1=v_\xi$  we will get nonzero eigenvectors
corresponding to all $\de_{\hw}$. Moreover,
$\om$ is nonzero at all of them. Let us check the latter. 
Indeed,  $x_{a_j}(\hu)=x_{\hu^{-1}(a_j)}(\xi)\neq t_j^{-1}$ 
in all intermediate 
$$\tPhi_j^{\hu} \ =\
T_j(x_{a_j}(\hu)-1) + (t_j^{1/2}-t_j^{-1/2}),$$
 since $\hu^{-1}(a_j)$ are positive.
However $\om$  vanishes after the application
of such $\tPhi_j^{\hu}$ only for such values.

We can always pick
a primitive character $\xi_o=\de_{\hu_o}$ satisfying
the same inequalities.  It  will be called {\it plus-primitive}
as well as the corresponding eigenvectors. 
The above argument gives more for the  $\xi=\xi_o$.
Since the eigenvector $1\in I_{\xi_o} $ is of multiplicity one in
$\H^o\ =\ I_{\xi_o}^\infty(\de_1)$, then
 $$
\tPhi_{\hw}(I_{\xi_o}^\infty(\de_1))\ =\
 I_{\xi_o}^\infty(\de_{\hw}).
$$ 
Really, the dimensions are the same and the image of $1$
is nonzero, since it belongs to 
any $X$-submodules of $I_{\xi_o}^\infty(\de_1)$. 
We note that this argument works for  any $\xi$  once
we know that the corresponding eigenvector  is simple.

Thus all eigenvectors of $I_{\xi_o}$ are exactly the
$\tPhi$-images of $1$ (in particular, they are simple),
and $\om$ is nonzero at them.
If $I_{\xi_o}$
contains a submodule where $\om$ vanishes, than the latter possesses
at least one eigenvector. It is impossible and
 $I_{\xi_o}$ is co-spherical. 

%\vfil%%% 
To go from $I_{\xi_o}$ to $I_\xi$, 
we need the following general lemma, where all the modules
are from the category $\o$ (we will apply it to  subquotients of 
induced representations).

%\vfil%%%
\proclaim{ Lemma}
Any submodule $V$ of co-spherical module $I$ is co-spherical.
If here $V\neq\{0\}$, then $I/V$ is not co-spherical. There 
exists a unique irreducible nonzero submodule  $U\subset I$.
A module with at least one  $\h_Y$-invariant homomorphism $\om\neq 0$
posseses a unique nonzero co-spherical quotient.
\label\COSPH\theoremnum*
\endproclaim

%\vfil%%%
{\it Proof.} 
The first and the second claims readily follow from the definition.
If there are two irreducible submodules $U,U'\subset I$,
then co-spherical $U$ is contained in the  co-spherical $\hat{U}=
U\oplus U'\subset I$.
Hence $\hat{U}/U\simeq U'$ couldn't be co-spherical (cf. Lemma 2.7
from [C5]). The last statement is  obvious as well. 
The kernel of the homomorphism to any co-spherical
module contains (the sum of) all submodules 
belonging to $\hbox{\ Ker}(\om)$. It cannot be bigger because
of the second assertion.

\proofbox

Let us check that $I_\xi$ is co-spherical. There exists a map
$I_\xi\to I_{\xi_o}$ sending $v_\xi =1$ to an eigenvector
of $I_{\xi_o}$ with $\xi$ as the eigenvalue. 
Here and further, all maps will be $\HH$-homomorphisms.
Since $I_{\xi_o}$
is co-spherical, the kernel
of this map belongs to $\hbox{\ Ker}(\om)$. For every $\de_{\hw},$
 there exist at least one 
eigenvector in $I_\xi$ apart from  $\hbox{\ Ker}(\om)$ (we have
already established this). Hence the image of $I_\xi$ contains
all eigenvectors of  $I_{\xi_o}$ and is surjective. Since the
spaces $I^\infty(\mu)$ for fixed $\mu$ have the same 
dimensions in all induced
modules with the characters from the same orbit, this map has
to be an isomorphism. It also proves that all co-spherical
modules (from the same orbit) are isomorphic.

\vfil%%%
 As to c), the lemma gives the uniqueness of the irreducible
submodule $U_\xi$ and that there are no biger modules (all of
them are co-spherical) with 
co-spherical quotients.
Any co-spherical irreducible
module $U$  containing an
 eigenvector with $\de_{\hu}$ as the character is the image
of a surjective homomorphism $I_{\hu}\to U$.
We again use the universality of the induced representations.
On the other hand, we can map $I_{\hu}$ to 
$I_\xi$ (since the latter contains an eigenvector 
corresponding to $\de_{\hu}$). The image will be co-spherical.
Therefore this map goes through $U$ (the last claim of the lemma)
and $U$ is isomorphic to $U_\xi$.
\proofbox
%\vfil%%%

\proclaim{ Corollary}
Imposing conditions (\ref\neqone), 
an arbitrary module $\De_\xi$ is isomorphic  to the module  $I_{\xi_o}$ 
for plus-primitive $\xi_o$.
In particular, it is generated
by any eigenvector with the character satisfying  (\ref\indcosph).
For instance, $\De(-\rho)=\De_{t^{-\rho}}$ (Proposition \ref\DELTAO)
is generated by $\de_{w_0}= t^{\rho}$ for the longest element
$w_0\in W$ and $\De_\#$ is the unique irreducible submodule
of $\De(-\rho)$, provided that $t$ is generic. 
\label\DESPH\theoremnum*
\endproclaim   
{\it Proof.}
There exists a nonzero homomorphism $I_{\xi_o}\to \De_\xi$.
Since both are co-spherical it has to be an isomorphism.

\proofbox

 {\bf Invariant forms}.
There are examples when co-spherical irreducibles exist
but there are no co-spherical induced representations at
all (say, for negative integers $k$). So the theorem does
not cover all of them (even if the orbit contains
a primitive character, which is always imposed). 
The following theorem makes the 
picture more complete and will be used to endow 
irreducible co-sperical representations with
$\ast$-invariant forms. 

\proclaim{ Theorem}
a) Let us assume that there exists a character $\xi'\in W^b(\xi)$ 
satisfying the condition dual to 
(\ref\indcosph):
$$
\eqalignno{ 
&x_{\ta}(\xi')\neq t_{\al} \for \hbox{ all\ }
\tal\in R_+^a.
&\eqnu
\label\indsph\eqnum*
}
$$ 
Then $I_{\xi'}$ has a unique 
nonzero irreducible quotient $U_{\xi'}$. It is co-spherical.
Other irreducible constituents (subquotients) 
are not co-spherical and do not contain 
eigenvectors corresponding to $\xi'$.

b) All  $I_{\tilde{\xi'}}$  are isomorphic to
 $I_{\xi'}$ 
for the characters $\tilde{\xi'}$ satisfying (\ref\indsph)
and  $W^b$-conjugated to $\xi'$. 
Otherwise $I_{\tilde{\xi}}$ has an irreducible quotient that
is not co-spherical.
The eigenvectors of $I_{\xi'}$ 
are simple. An eigenvector 
generates  $I_{\xi'}$ (and $U_{\xi'}$) if and only if its
character belongs to $\{\tilde{\xi'}\}$.   
Any co-sperical irreducible representation $U$ with a character
from the orbit of $\xi'$ 
is isomorphic to $U_{\xi'}$. 
\label\IRRCOSPH\theoremnum*
\endproclaim

{\it Proof.} Primitive characters satisfying 
(\ref\indsph) will be refered to as {\it minus-primitive}.
They always exist in the orbit of the character $\xi'$
(satisfying (\ref\indsph) and conjugated to a primitive one). 
The same terminology will be used for the corresponding
eigenvectors.
Actually $I_{\xi'}$ are  counterparts of {\it spherical}
induced representations from the theory of affine Hecke algebras.
So we can follow  [C5] closely. 

Let $v$ be 
an $X$-eigenvector of an irreducible co-spherical
$U$ with the character $\de_{\hu}$. 
Starting with $v$, we can construct the eigenvector 
$v_1=\tPhi_{j_1}^{\hu}(v)$ with the character $\de_{\hu_1}$
for $\hu_1=s_{j_1}\hu$, then  $v_2=\tPhi_{j_2}^{\hu_1}(v_1)$ and so on,
till we get minus-primitive $v'=v_l$
associated with some $\xi'_o$.
As in Lemma 2.13 from [C5], the values
$$
x_{j_1}(\hu),\ x_{j_2}(\hu_1),\ldots,\ x_{j_l}(\hu_{l-1})
$$
can be chosen avoiding $t_{j_r}^{-1} \ (1\le r\le l)$.

Let us check that $v'\neq 0$. Otherwise $v_r\neq 0$ in this chain
for some $r$ and 
$$
0=v_{r+1}\ =\  
(t_{j_r}-1)(T_{j_r} +t_{j_r}^{-1/2})(v_{r}),\where
t_{j_r}\neq -1.
$$
However $U=\h_Y(v_{r})$, since $v_r$ is an eigenvector 
and $U$ is irreducible. Hence,  $\om(v_r)=0$  and
$\om(U)=0$, which is impossible. 

This gives a  surjective map $I_{\xi'_o}\to U$. 
It  establishes an isomorphism between
$\H^o= I_{\xi'_o}^\infty(\de_1)$ and $U^\infty(\xi'_o)$
because the first is an irreducible $\h_X^o$-module
(or since the eigenvector $1$ is simple in $\H^o$).
Thus the kernel of this map (and any its subquotients)
cannot contain $\xi_o'$-eigenvectors.
Hence $U$ is a unique irreducible quotient of $I_{\xi'_o}$.

By the way, here the combinatorial part can be simplified
a bit. We can finish with $v'$ associated with any
primitive character. Indeed (cf. the proof of Theorem
\ref\INDCOSPH, b)), it has to satisfy (\ref\indsph)
because otherwise we can construct a quotient of $I_{\xi'_o}$
which contains no $\om\neq 0$ (use a proper intertwiner).

Vice versa, if $\xi'_o$ is minus-primitive then $ I_{\xi'_o}$ has a
unique irreducible quotient $U_{\xi'_o}$
(it holds true for all primitive
characters). On the other hand, it has a unique
(nonzero) co-spherical
quotient $V$, which contains a co-spherical irreducible
submodule $U$. It is a quotient of a proper $ I_{\tilde{\xi'_o}}$.
Due to [C5], Lemma 2.8, the latter is isomorphic to $ I_{\xi'_o}$
(any  primitive eigenvectors   are
connected by invertible intertwiners).
So we get a nonzero homomorphism from $ I_{\xi'_o}$
onto $U$, which results in  $V=U$ and $U\simeq U_{\xi'_o}$. 

To establish the necessary isomorphism without any reference to
[C5] and  check the remaining statements of the theorem
it is convenient to  apply Theorem \ref\INDCOSPH
for the character of $\h_Y$ sending
$$
T_i\to -t_i^{-1/2},\ \pi_r\to 1,\  0\le i\le n,\  r\in O.
$$
The modules under consideration become co-spherical
for this character. 

First, it gives that all eigenvectors of $I_{\xi'_o}$
are simple and can be obtained from $1$ by the 
intertwiners. However now once the intertwiner is not invertible
its image  belongs to $\hbox{\ Ker}(\om)$. 
Second, all $I_{\xi'}$ satisfying (\ref\indsph) are isomorphic
to $I_{\xi'_o}$ (and to each other). Moreover,
all eigenvectors in $I_{\xi'_o}$ 
corresponding to $\xi'$ are connected with $1$ by invertible
intertwiners (and only them).
Then if an irreducible constituent of $I_{\xi'_o}$ 
contains an eigenvector associated with $\xi'$ then 
 $I_{\xi'_o}$ maps through $I_{\xi'}$ onto it. Hence it can happen
for  $U_{\xi'_o}$ only. 

\proofbox

Let us assume that both plus and minus-primitive characters
$\xi_o,\xi_o'$ belong to the orbit of $\xi$.
Then $U_\xi\simeq U_{\xi'}$ (since they are unique co-sperical
irreducible constituents) and moreover the first
contains the whole $I^\infty_\xi(\xi')$.
The modules $I_{\xi_0}, \ I_{\xi_o'}$ are dual to each other
 in the following sense.
Given a module $I=\oplus I^\infty(\mu)\in \o$, we 
combine the anti-involution $\diamond$ from
(\ref\antivee) with the natural anti-action
of $\HH\ $ on the
$$
\eqalign{
& I^{\diamond}\equal \{f\in \hbox{Hom}(I,\C_{q,t}) \hbox{
\ s.t.\ } f(I^\infty(\mu))=0\hbox{\ for\ almost\ all\ } \mu
\}.
}
\eqnu
\label\ocatdual\eqnum*
$$
We claim that $I_{\xi_o}^\diamond\simeq I_{\xi'_o}$.

First of all, $U_{\xi_o}^\diamond\simeq U_{\xi'_o}$.
Indeed, they have the same set of characters and coinciding
dimensions of the spaces of generalized vectors 
(it is true for any dual modules). Hence  $U_{\xi_o}^\diamond$
can be covered by $I_{\xi'_o}$. However the latter has a unique
irreducible quotient which is just $U_{\xi'_o}$.
The map $\nu:I_{\xi'_o}\to  I_{\xi_o}^\diamond$
sending $1$ to an eigenvector corresponding to $\xi'_o$
composing with the map $ I_{\xi_o}^\diamond \to 
U_{\xi_o}^\diamond$ is obviously nonzero.  
The latter map is a dualization of the embedding
$U_{\xi_o}\subset  I_{\xi_o}$ (Theorem \ref\INDCOSPH). 
Similarly, the  
module  $ I_{\xi_o}^\diamond$
has a unique nonzero irreducible quotient.
Hence $\nu$ is surjective and has to be an isomorphism because
the dimensions of the spaces $I^\infty(\mu)$ are
the same for $I_{\xi'_o}$ and   $I_{\xi_o}$ (their characters
are from the same orbit). 

\proclaim{ Corollary}
In the setup of Theorem \ref\IRRCOSPH,
the module $I_{\xi'}$ 
possesses a
unique (up to proportionality) $\ast$-invariant form 
in the sense of Proposition \ref\YONE,a).
Its radical $\r$ is exactly the kernel $\k$ of the map 
$I_{\xi'}\to U_{\xi'}$ from Theorem \ref\IRRCOSPH. The restriction
is a unique (nondegenerate) $\ast$-invariant form on $U_{\xi'}$.
\endproclaim
\label\FORMSPH\theoremnum*
{\it Proof.}
Any $I_{\xi'}$ can be considered as a limit of
a one-parametric family 
of proper $F_{\tilde{\xi'}}$ ensuring the same inequalities
$\tilde{\xi'}$. 
Therefore $I_{\xi'}$  has a nonzero $\ast$-invariant form.
If $\r$ is less than $\k$, then $I_{\xi'}/\r$ contains an
irreducible submodule $V\not\simeq U_{\xi'}$. 
All generalized $X$-eigenvectors
with the characters not from $V$ belong to its
orthogonal compliment $V'$. Since $V$ has no eigenvectors
associated with $\xi'$,  the image of $1$ belongs to $V'$.
A contradiction.
The uniquiness of the $\ast$-invariant form on $U_{\xi'}$
follows from the irreducibility.
\proofbox 

Provided (\ref\neqone), we claim that
$I_{\xi'_o}\simeq F_\xi$ for minus-primitive $\xi'_o$
from the orbit of $\xi$ (cf. Corollary \ref\DESPH).
Indeed, one can map $I_{\xi'_o}$ into $F_\xi$ and 
replace $\De_\xi$ by $I_{\xi_o}$ for plus-primitive 
$\xi_o$. Then 
the pairing $\{\ ,\ \}$ can be extended
to $I_{\xi'_o}\times I_{\xi_o}$. 
The right radical ($\subset I_{\xi_o}$) has to contain $U_{\xi_o}$
(the smallest irreducible) if the resulting pairing is degenerate. 
However it is imposible because the image of $I_{\xi'_o}$ in $F_\xi$
contains $U_{\xi'_o}\simeq U_{\xi_o}$ as a constituent.  
We see that the  above corollary generalizes the
calculation of the radical of the form $\lan\ ,\ \ran_1$ on $F_\xi$ for
$\xi=t^{-\rho}$.

We also note that 
Corollary \ref\ISOVD ,b) (which is almost equivalent to the Main
Theorem) readily follows from the theory of co-spherical
representations.
Applying $\vep$ to $\v_0$ one gets an irreducible and 
$X$-co-spherical representation (it possesses $\om\neq 0$)
with $t^{-\rho}$ as an eigenvalue.
Hence it is isomorphic to $U_{t^\rho}$,
which in its turn  is isomorphic to $\De_{\#}$ (for generic $t$).

\proclaim{ Proposition}
a) Given a finite dimensional irreducible $\h_X$-module $U$ let
as assume that
$$
\eqalignno{ 
&1\neq q_a^j x_{\ta}(\xi)\neq t_{\al}^{\pm 1} \for \hbox{ all\ }
\tal\in R^a,\ \Z\ni j>0,
&\eqnu
\label\indirU\eqnum*
}
$$
where $\xi$ is any eigenvalue of $U$ (does not matter which because
they are $W$-conjugated). Then the induced $\HH\ $-module 
$M_U=\hbox{Ind\,}_{\h_X}^{\HH}(U)$ is irreducible.

b) % $q$ is an indeterminate then 
Any irreducible $\HH\ $-module $M$
posessing an $X$-eigenvector with the eigenvalue $W^b$-conjugated to
$\xi$ from (\ref\indirU) is induced 
from its irreducible  $\h_X$-submodule $U$ 
generated by an eigenvector with the character from
$W(\xi)$. Such a submodule $U$ is unique 
(if the orbit $W(\xi)$
is fixed),  $M_U\simeq M_{U'} \Rightarrow U\simeq U'$.
\label\INDIRU\theoremnum*
\endproclaim
{\it Proof}
 is close to the  proof of the irreducibility from 
Theorem \ref\INDCOSPH. The left hand side inequality gives that
$M^\infty(\xi)=U^\infty(\xi)$ for $M=M_U$, where 
(see (\ref\xeigen)) by $M^\infty(\xi)$ we mean the space of all generalized
eigenvectors in $M$ associated with a character $\xi$. 
Indeed, it is true for $I_\xi$ and for the induced $\h_X$-module
$I_\xi^X=\Indcx^{\h_X}(\xi)$ instead of $M$ and $U$. 
However $I_\xi$ and $I_\xi^X$
cover $M,U$ naturally (the same holds for the spaces of generalized
eigenvectors). Once the coincidence is true at $\infty$-level it
is valid for all levels. In particular,  $M^1(w(\xi))=U^1(w(\xi))$
for any $w\in W$. This argument also gives that 
$$
\eqalignno{ 
&U\ =\ \oplus_{\ze\in W(\xi)} M_U^\infty(\ze).
&\eqnu
\label\indUsum\eqnum*
}
$$

Given $b\in B$ and  $w\in W$,
$$
\eqalignno{ 
&(\tPhi_{\pi_b}^{\ze})^{-1}(M^1(\pi_b(\ze))\ =\ M^1(\ze),\cr
&\where \tPhi_{\pi_b}^{\ze}=
\tPhi_{\pi_b}(\ze),\ \ze\in W(\xi).
&\eqnu
\label\indUeig\eqnum*
}
$$
We follow the notations from  
(\ref\Phiind) and use the invertibility of $ \tPhi_{\pi_b}^{w(\xi)}$
in the space $M^1(w(\xi))$ (or even in $M^\infty(w(\xi))$) thanks to
the right hand side inequality. The elements $\{\pi_b(\ze)\}$ 
constitute the whole orbit $W^b(\xi)$.
 Any nonzero irreducible submodule of $M$
has at least one $X$-eigenvector for a character from
the orbit $W^b(\xi)$. Due to (\ref\indUeig) it generates the whole
$M$. 
 
An arbitrary irreducible
$M$ posessing an eigenvalue $W^b$-conjugated to $\xi$ from
 (\ref\indirU) can be represented as $M_U$ where 
$U$ is any irreducible
$\h_X$-submodule of $M$ with an eigenvalue from $W(\xi)$.
The existence of $U$ results from  the same formula 
(\ref\indUeig).
Moreover $U$ can be reconstructed uniquely
as a submodule of $M$ by means of (\ref\indUsum).
Therefore  $M\simeq M' \Rightarrow U\simeq U'$.
\proofbox

Let us discuss the structure of $M_U$ upon the restriction
to $\h_X$, provided (\ref\indirU).
 First of all, $M_U=\oplus J_{c,\ze}$ where
$c\in B_-$, $\ze$ runs over a fixed  set of representatives
of $W(\xi)\hbox {\,mod\,} W_c$ for the centralizer $W_c$ of $c$ in $W$,
$$
\eqalignno{ 
&J_{c,\ze}\ =\ \oplus_{\txi} M_U^\infty(\txi),\
\txi\in W(c(\ze)).
&\eqnu
\label\jcze\eqnum*
}
$$
All $\{J\}$ are $\h_X$-submodules. Their structure can be described
as follows:
$$
\eqalignno{ 
&J_{c,\ze}\ =\ \oplus_{\om}\tPhi_\om (\tU_c),\  
 \tU_c\ =\ \tPhi_{c}(U_c),\cr
&U_c\equal\oplus_{\xi'} U^\infty(\xi')),
\where \xi' \in W_c(\ze), 
&\eqnu
\label\jcze\eqnum*
}
$$
$$
\{\om=\om_b^{-1},\ b\in W(c)\} = \{\om\in W,\ \la(\om)\subset
\{\al\in R_+ \hbox{\ s.t.\ } (\al,c)<0\}\}.   
$$
Here  the intertwiners $\tPhi_\om, \tPhi_c$ are invertible
because $c=c_-=\pi_c$ and $(\al,c)\neq 0$ for all $\al\in \la(\om)$.
The space $U_c$ is a module over the subulgebra $\p_c$
of $\h_X$ generated by $\{T_i, s_i(c)=c\}$ and $\C[X]$.
So does $\tU_c$. Indeed, to  apply $\Phi_c$ means to replace the
action of $X_j$ by that of $q^{(c,b_j)}X_j$ for all
$1\le j\le n,\ s_j(c)\neq c$ without changing $\{T_i\}$.

Finally, $J_{c,\ze}$ is isomorphic to the representation
of $\h_X$ induced from the $\p_c$-module $\tU_c$.
Its irreducibllity is equvalent to the irreducibility of
the $\p_c$-module $U_c$. The simplest example is the
decomposition of $\v_0$
as an $\h_Y$-module. It was considered in [C4] (formulas (3.15)-(3.17)).   

Summarizing, in the case of generic $q$ the classification of
irreducible representations of $\HH\ $ is not far from 
that in the affine case. If all $t_\nu$ coincide then we can
use directly the main theorem from [KL] and moreover try
to generalize it to arbitrary $q$. The latter seems to be quite
possible because $\HH\ $ has a natural $K$-theoretic interpretation
due to [KK] and more recent [GH],[GKV]. The list of
finite groups which are
expected to appear in the data (see [KL]) 
can be rather complicated.

\appendix{}{Degenerate double affine Hecke algebras}
The theory of induced and co-spherical representations
is very close to  that of the degenerate ones.   
On the other hand, the degenerate case
is not self-dual, which makes quite different induced 
(basic or functional) representations
associated with  $\C[x]$ and $\C[y]$. 
To connect them (as we did many times 
in the paper) one needs to go to the difference theory.

Let us fix 
$\ka=\{ \ka_\nu \in \C , \nu \in \nu_R,\ \eta\in \C\}$  and
introduce a linear function $\rho_{\ka}$ on $[a,u]\in 
\R^{n}\times \R$
setting
 $\rho_{\ka}(a_j)= \ka_j$
for  $0\le j\le n$. As always, $a_j=\al_j^\vee$. We will also
need a $\ka$-deformation of the Coxeter number:  $h_\ka=\ka_0+\rho_\ka(\th)$.

The {\it  degenerate (graded) double affine Hecke algebra} $\HH'_\eta$ is 
 algebraically generated by
 the group algebra $\C [W^b]$ and the pairwise commutative
$$
\eqalignno{
&y_{\ta}\equal
\sum ^n_{i=1}(a,\al_i)y_i - u\eta \for 
\ta=[a,u],
&\eqnu
}
$$ 
satisfying  the following relations:
$$
\eqalignno{
&s_j y_{\ta}-y_{\{s_j(\ta)\}}s_j\ =\ \ka_j(a,\al_j ),
\ 0\le j\le n, \cr
&\pi_r y_{\ta}\ =\ y_{\pi_r(\ta)}\pi_r, \ r\in O.
&\eqnu
\label\suka\eqnum*
}
$$

Without $s_0$ and $\pi_r$ we arrive at the defining relations
of the graded affine Hecke algebra from  [L] (see also [C8]).
It is a natural degeneration  
of the double affine
Hecke algebra when $q\to 1, t\to 1$ (see below).

 We will  use
the derivatives of $\C[x]$:
$
\ \partial_{[a,u]}(x_{[b,v]})=-(a ,b)x_{[b,v]}.  
$
Here the sign is minus to make the definition compatible
with  (\ref\saction). We note that 
$\tw(\partial _{\tb})=
\partial _{\tw(\tb)}, \  \tw\in W^b$.

\proclaim{ Theorem}
a) Taking $\tb=[b,v] \in \R^{n}\times \R$,
the following family
of the  operators
$$
\eqalignno{
\ty_{\tb} \equal
&\pa_b +
\sum_{\tal\in R^a_+} { \ka_{\al}(b,\al)\over
(X_{\tal^\vee}^{-1}-1) }
\bigl( 1-s_{\tal} \bigr)+ \rho_{\ka}(\tb)\ 
&\eqnu
\label\dunell\eqnum*
}
$$
is commutative and satisfies (\ref\suka) for $\eta=-h_\ka$.

b) Their non-affine counterparts
$$
\eqalignno{
y_{[b,v]} \equal
&\pa_b +
\sum_{\al\in R_+} { \ka_{\al}(b,\al)\over
(X_{\al^\vee}^{-1}-1) }
\bigl( 1-s_{\al} \bigr)+ \rho_{\ka}(b) - v,\ 
&\eqnu
\label\dunk\eqnum*
}
$$
are pairwise commutative and satisfy (\ref\suka) for $\eta=1$ and
the group
$$\tW^b\equal W\lsmash B_X, \where B_X\ =\ \{X_b,b\in B\}.
$$ 
The group $B_X$ acts naturally in $\C[x]$.
\label\DUNDEG\theoremnum*
\endproclaim
{\it Proof.} The first statement is from [C7], where the
convergence problem is managed in full detail
(the difference version also exists). Presumably it 
is a good starting point for the harmonic analysis in the Kac-Moody
case at critical level. The  second claim is essentially  from
[C8]. Since it can be easily deduced 
 from the difference theory, we will outline the proof. 

Setting 
$$q=1+h,\ t_j=q_j^{k_j},\ \ka_j=2k_j/\nu_j \ (\hbox{\ i.e.\ }
 t_j= q^{\ka_j}),\ Y_b=1+h y_b,
$$ let us tend $h$ to zero ignoring the terms of order $h^2$
(without touching $X_b,\ka_j$).
We will readily arrive at the relations
(\ref\suka) for $y_b$  with the constant
$\eta=1$. The formula for $y_b$ is exactly  (\ref\dunk).
\proofbox

The formulas for the intertwiners of $\HH'=\HH'_1$
 generalize those for the degenerate (graded) Hecke algebras
(see [L] and [C5,C8]) and result from the limiting
procedure.
One can  use them to create the Opdam and Jack polynomials
(the Macdonald ones in the differential setup).
Starting with $\Phi$ from (\ref\Phi), it is necessary  to apply
the involution $\vep$ and then  tend $t\to 0$:
$$
\eqalignno{
&\Phi'_i\ =\ 
s_i - {\ka_i \over y_{a_i} },\ 
\Phi'_0\ =\ 
X_\th s_\th + {\ka_0 \over y_{\th}+1 },\cr
&\pi_r'\ =\ X_r\om_r^{-1}, \for a_i=\al_i^\vee,\ 1\le i\le n,\ r\in O.
&\eqnu
\label\Phiprime\eqnum*
}
$$
The operators $\pi_r'$ in the case of $GL_n$ (they are of 
infinite order) play the key role in [KS].

Due to [C5], Corollary 2.5,  the intertwiners always  lead to
 Lus\-ztig's isomorphisms [L]. It gives that the
algebras $\HH$ and $\HH'$  are isomorphic after proper completion
for generic $q,t$ in the sense
of [L] (for formal parameters) or in the sense of [C5] (in
the category $\o$).   

The theory of the degenerate induced and 
co-spherical representations is very
close to what we did in the $q,t$-case. 
It is somewhat simplier because any 
characters are conjugated to primitive ones. 
The $y$-induced representation $J_\xi$
is generated by $v$ such that $y_b(v)=y_b(\xi)v$ 
for the character $\xi:\C[y_1,\cdots,y_n]\to \C$ in the above
notations. It is
co-spherical if and only if
$$
\eqalignno{
& -\ka_\al-1-\Z_+\not\ni y_a(\xi)\not\in \ka_\al+\Z_+,\
a=\al^\vee\in R_+^\vee.
&\eqnu
\label\Phiprime\eqnum*
}
$$
As for  irreduciblity,  the inequalities have to 
hold for $\pm \ka_\al$ instead of $\ka_\al$. The proof remains the same. 

We note that this paper mainly follows the same lines
as  the $p$-adic theory (although
many properties of the double affine Hecke algebras are brand
new). Probably the most interesting point is that our methods
(essentially $p$-adic) work very well in 
the differential theory via the semi-classical limit.
Since the latter is a  generalization of the classical harmonic 
analysis in the zonal case, 
we have a new foundation for
the Harish-Chandra theory of spherical functions.

%
%
%
%      REFERENCES
% 
%
%
%\vskip 15pt
\AuthorRefNames [BGG]
\references
%\medskip
%\ninerm
%\baselineskip=10pt %!
\vfil

[A]
\name {K. Aomoto},
{On product formula for Jackson integrals associated with root
systems}, Preprint (1994).

[BGG]
\name {I.N. Bernstein}, \name{I.M. Gelfand}, and
\name{S.I.  Gelfand},
{Schubert cells and the cohomology of $G/P$}, Russ.Math.Surv.
{28} (1973) ,1--26.

[B]
\name{N. Bourbaki},
{ Groupes et alg\`ebres de Lie}, Ch. {\bf 4--6},
Hermann, Paris (1969).

[C1]
\name{I. Cherednik},
{ Double affine Hecke algebras, 
Knizhnik- Za\-mo\-lod\-chi\-kov equa\-tions, and Mac\-do\-nald's 
ope\-ra\-tors},
IMRN (Duke M.J.) {  9} (1992), 171--180.

[C2]
\bibline,{ Double affine Hecke algebras and  Macdonald's
conjectures},
Annals of Mathematics {141} (1995), 191--216.

[C3]
\bibline, 
{ Macdonald's evaluation conjectures and
difference Fourier transform},
Inventiones Math. {122} (1995),119--145.

[C4]
\bibline, 
{ Nonsymmetric Macdonald polynomials },
IMRN {10} (1995), 483--515.

[C5]
\bibline, 
{A unification of Knizhnik--Zamolodchikov
and Dunkl operators via affine Hecke algebras},
Inventiones Math. {  106}:2  (1991), 411--432.

[C6]
\bibline,
{Induced represenations of double affine Hecke algebras and
applications}, Math. Res. Let. {1} (1994), 319--337. 

[C7]
\bibline,
{Elliptic quantum many-body problem and double affine
 Knizhnik - Zamolodchikov equation},
Commun. Math. Phys. {169}:2  (1995), 441--461.

[C8]
\bibline,
{ Integration of quantum many- body problems by affine
Knizhnik--Za\-mo\-lod\-chi\-kov equations}, 
Advances in Math. {106}:1 (1994), 65--95
(Pre\-print RIMS--{  776} (1991)).

[GG]
\name{H. Garland}, and \name {I. Grojnowski}, 
{Affine Hecke algebras associated to Kac-Moody groups},
Preprint q-alg/9508019 (1995).

[GKV]
\name{V. Ginzburg} , \name{M. Kapranov} , 
and \name{E. Vasserot},
{Residue construction of Hecke algebras},
Preprint alg-geom/9512017 (1995).

[HO]
\name{G.J. Heckman}, and \name{E.M. Opdam},
{Harmonic analysis for affine Hecke algebras},
Preprint (1996).

[IM]
\name{N. Iwahori}, and \name{H. Matsumoto},
{On some Bruhat decomposition and the structure of
Hecke rings of $p$-adic Chevalley groups},
Publ.Math.  IHES {25} (1965), 5--48.

[Ito]
\name{M. Ito}, {On a theta product formula
for Jackson integrals associated with root
systems of rank two},
Preprint (1996).

[K]
\name {V.G. Kac},
{Infinite dimensional Lie algebras},
Cambridge University Press, Cambridge (1990).

[Kn]
\name {F. Knop},
{Integrality of two variable Kostka functions},
Preprint (1996). 

[KS]
\name {F. Knop}, and \name{S. Sahi},
{A recursion  and a combinatorial formula for Jack
polynomials},
Preprint (1996).

[KL]
\name{D. Kazhdan}, and \name{ G. Lusztig},
{  Proof of the Deligne-Langlands conjecture for Hecke algebras},
Invent.Math. {  87}(1987), 153--215.

[KK]
\name{B. Kostant}, and \name{ S. Kumar},
{  T-Equivariant K-theory of generalized flag varieties,}
J. Diff. Geometry{  32}(1990), 549--603.

[LV]
\name{L. Lapointe}, and \name{L. Vinet},
{A Rodrigues formula for the Jack polynomials and the
Macdonald-Stanley Conjecture},
IMRN {9} (1995), 419--424.

[L]
\name {G. Lusztig},{ Affine Hecke algebras and their graded version},
J. of the AMS { 2}:3 (1989), 599--685.

[M1]
\name{I.G. Macdonald}, {  A new class of symmetric functions },
Publ.I.R.M.A., Strasbourg, Actes 20-e Seminaire Lotharingen,
(1988), 131--171 .

[M2]
\bibline, {  Orthogonal polynomials associated with root 
systems},Preprint(1988).

[M3]
\bibline, { Affine Hecke algebras and orthogonal polynomials},
S\'eminaire Bourbaki{  47}:797 (1995), 01--18.

[M4]
\bibline, { A formal identity for affine root systems},
Preprint (1996).

[O1]
\name{E.M. Opdam}, 
{  Some applications of hypergeometric shift
operators}, Invent.Math.{  98} (1989), 1--18.

[O2]
\bibline, {Harmonic analysis for certain representations of
graded Hecke algebras}, 
Preprint Math. Inst. Univ. Leiden W93-18 (1993).

[R]
\name{J. Rogawski},
{On modules over the Hecke algebra of a $p$-adic group},
Invent. Math. {79} (1985), 443--465.

[S] 
\name {S. Sahi},
{Interpolation, integrality, and a generalization of
Macdonald's polynomials},
Preprint (1996).

[StK]
\name {J.V. Stokman}, and \name {T.H. Koornwinder},
{Limit transition for $BC$ type multivariable orthogonal
polynomials}, Canad. J. Math. (to appear).

\endreferences
\bye